\documentclass[12pt,letterpaper,reqno]{amsart}
\usepackage{amsfonts}
\usepackage{color}
\usepackage{epsfig}
\usepackage{fancyhdr}
\usepackage[latin1]{inputenc}
\usepackage{amsmath}
 \usepackage{mathrsfs}
\usepackage{thmdefs}
\usepackage{natbib}
\usepackage{graphicx}
\usepackage[colorlinks=true,linkcolor=red,citecolor=blue]{hyperref}
\usepackage{amssymb}
\usepackage{dsfont}
\usepackage{fancyhdr}
\usepackage{graphicx}

\theoremstyle{definition} \theoremstyle{remark}
\numberwithin{equation}{section}

\renewcommand{\cite}{\citet}
\numberwithin{equation}{section}

\newtheorem{assumption}{Assumption}

\begin{document}

\author{Michel BRONIATOWSKI$^{1}$ and Amor KEZIOU$^{1,2}$}
\address{$^{1}$LSTA, Universit\'e Pierre et Marie Curie -
Paris 6. E-Mail: michel.broniatowski@upmc.fr
 \newline
 $^{2}$Laboratoire de Math\'ematiques de Reims, EA 4535, Universit\'e de Reims Champagne-Ardenne.
  E-Mail: amor.keziou@upmc.fr}
\date{March 2011}
\title[Estimation and Test Under Moment Condition
Models]{Divergences and Duality for Estimation and Test under
Moment Condition Models}

\maketitle

\begin{abstract}
We introduce estimation and test procedures through divergence
minimization for models satisfying linear constraints with unknown
parameter. These procedures extend the empirical likelihood (EL)
method and share common features with generalized empirical
likelihood  approach. We treat the problems of existence and
characterization of the divergence projections of probability
distributions on sets of signed finite measures. We give a precise
characterization of duality, for the proposed class of estimates
and test statistics, which is used to derive their  limiting
distributions (including the EL estimate and the EL ratio
statistic) both under the null hypotheses and under alternatives
or misspecification. An approximation to the power function is
deduced as well as the sample size which ensures a desired power
for a given alternative.
 \vspace{2mm} 
\end{abstract}


\tableofcontents

\section{Introduction and notation}

\noindent Statistical models are often defined through estimating
equations
\begin{equation*}
\mathbb{E} \left[g(X,\theta)\right]=0,
\end{equation*}
where $\mathbb{E} [\cdot]$ denotes the mathematical expectation,
$g:=(g_1,\ldots,g_l)^{\top}\in\mathbb{R}^l$ is some specified
vector valued function of a random vector $X\in \mathbb{R}^m$ and
a parameter vector $\theta\in\Theta\subset \mathbb{R}^d$. Examples
of such models are numerous, see e.g.
\noindent\cite{Qin-Lawless1994}, \cite{Haberman1984},
\cite{Sheehy1987}, \cite{McCullagh_Nelder1983}, \cite{Owen2001}
and the references therein. Denoting $M^{1}$ the collection of all
probability measures (p.m.) on the measurable space
$(\mathbb{R}^m,\mathcal{B}(\mathbb{R}^m))$, the submodel
$\mathcal{M}_{\theta}^{1}$, associated to a given value $\theta$
of the parameter, consists of all distributions $Q$ satisfying $l$
linear constraints induced by the vector valued function
$g(.,\theta)$, namely
\begin{equation*}
\mathcal{M}_{\theta}^{1}:=\left\{  Q\in M^{1}~\text{ such that
}~\int g(x,\theta)~dQ(x)=0\right\},
\end{equation*}
with $l\geq d$. The statistical model which we consider can be
written  as
\begin{equation}\label{modele}
\mathcal{M}^{1}:=\bigcup_{\theta\in\Theta}\mathcal{M}_{\theta}^{1}.
\end{equation}

\noindent Let $X_{1},...,X_{n}$ denote an i.i.d  sample of $X$
with unknown distribution $P_{0}$. We denote $\theta_{0}$, if it
exists, the value of the parameter such that $P_{0}$ belongs to
$\mathcal{M}_{\theta_{0}}^1$, namely the value satisfying
 $$\mathbb{E}\left[g(X,\theta_0)\right]=0,$$ and we assume obviously that $\theta_0$
 is unique. This paper addresses the two following natural
questions:

\textit{Problem 1}: Does $P_{0}$ belong to the model
$\mathcal{M}^{1}$?

\textit{Problem 2}: When $P_{0}$ is in the model, which is the
value $\theta_{0}$ of the parameter for which
$\mathbb{E}\left[g(X,\theta_{0})\right]=0$? Also can we perform
tests about $\theta_{0}$? Can we construct confidence areas for
$\theta_{0}$?\\

\noindent We note that these problems have been investigated by
many authors. \cite{Hansen1982} considered generalized method of
moments (GMM). \cite{Hansen_Healton_Yaron1996} introduced the
continuous updating (CU) estimate.   The empirical likelihood (EL)
approach, developed by \cite{Owen1988} and \cite{Owen1990}, has
been investigated in the context of model (\ref{modele}) by
\cite{Qin-Lawless1994} and \cite{Imbens1997} introducing the EL
estimate. The recent literature in econometrics focusses on such
models; \cite{Smith1997}, \cite{NeweySmith2004} provided a class
of estimates called generalized empirical likelihood (GEL)
estimates which contains the EL and the CU ones.
\cite{Schennach2007} discussed the asymptotic properties of the
empirical likelihood estimate under misspecification; the author
showed the important fact that the EL estimate may cease to be
root $n$ consistent when the functions $g_j$ defining the moments
conditions and the support of $P_0$ are unbounded. Among other
results pertaining to EL, \cite{NeweySmith2004} stated that EL
estimate enjoys optimality properties in term of efficiency when
bias corrected among all GEL estimates including the GMM one.
Moreover, \cite{Corcoran1998} and \cite{Baggerly1998} proved that
in a class of minimum discrepancy statistics (called power
divergence statistics), EL ratio is the only one that is Bartlett
correctable. Confidence areas for the parameter $\theta_{0}$ have
been considered in the seminal paper by \cite{Owen1990}. Problems
1 and 2 have been handled via EL and GEL approaches  in
\cite{Qin-Lawless1994}, \cite{Smith1997} and \cite{NeweySmith2004}
under the null hypothesis $\mathcal{H}_0 : P_0\in\mathcal{M}^1$;
the limiting distributions of the GEL estimates and the GEL test
statistics have been obtained under the model and under the null
hypotheses. \cite{Imbens1997} discusses the asymptotic properties
of the EL and exponential tilting estimates under misspecification
and give the formula of the asymptotic variance, using dual
characterizations,  without presenting the hypotheses under which
their results hold. \cite{Chen_Hong_Shum2007} give the limiting
distribution of the EL estimate under misspecification as well as
the EL ratio statistic between a parametric model and a moment
condition model. The paper by \cite{Kitamura2007} gives a
discussion of duality for GEL estimates under moment condition
models. \cite{Bertail2006} uses duality to study, under the model,
the asymptotic properties  of the EL ratio statistic and its
Bartlett correctability; the author extends his results to
semiparametric problems with
infinite-dimensional parameters.\\

\noindent The main contribution of the present paper is the
precise characterization of duality for a large class of estimates
and test statistics (including GEL and EL ones) and its use in
deriving the limiting properties of both the estimates and the
test statistics under misspecification and under  alternatives
hypotheses. Moreover,
\begin{enumerate}
\item [1)] The approach which we develop is based on
\textit{minimum discrepancy estimates}, which extends the EL
method and has common features with minimum distance and GEL
techniques, using merely divergences. We present a wide class of
estimates,  test statistics and confidence regions for the
parameter $\theta_{0}$ as well as various test statistics for
\textit{Problems} 1 and 2,
all depending on the choice of the divergence.\\

\item [2)] The limiting distribution of the EL test statistic
under the alternative and under misspecification remains up to
date an open problem. The present paper fills this gap; indeed, we
give the limiting distributions of the proposed estimates and test
statistics (including the EL ones) both under the null hypotheses,
under alternatives and under
misspecification.\\

\item [3)] The limiting distributions of the test statistics under
the alternatives and misspecification are used to give an
approximation to the power function and the sample size which
ensures a desired power for a given alternative.\\

\item [4)] We extend confidence region (C.R.)  estimation
techniques based on EL (see \cite{Owen1990}), providing a wide
range of such C.R.'s, each one depending upon a specific
divergence.
\end{enumerate}

\noindent From the point of view of the statistical criterion
under consideration, the main advantage, of using a divergence
based approach and duality, lays in the fact that it leads to
asymptotic properties of the estimates and test statistics under
the alternative, including misspecification, which cannot be
achieved through the classical EL context. In the case of
parametric models of densities, \cite{White1982} studied the
asymptotic properties of the parametric maximum likelihood
estimate and the parametric likelihood ratio statistic under
misspecification; \cite{Keziou2003} and \cite{Broniatowski_Keziou2009} stated the consistency and obtained the limiting distributions of the minimum
divergence estimates and the corresponding test statistics
(including the parametric likelihood ones) both under the null
hypotheses and the alternatives, from which they deduced an
approximation to the power function. In this paper, we extend the
above results for the proposed class of estimates and test
statistics (including the EL ones) in the context of
semiparametric models
(\ref{modele}). \\

\noindent The rest of the paper is organized as follows. Section 2
describes the statistical divergences used in the sequel. Section
3 is devoted to the description of the proposed estimation and
test procedures. In Section 3, we adapt the Lagrangian duality
formalism  to the context of statistical divergence, and we use it
to give practical formulas (for the study and the numerical
computation) of the proposed estimates and test statistics.
Section 5 deals with the asymptotic properties of the estimates
and the test statistics under the model and under
misspecification. Simulations results are given in Section 6. All
proofs are postponed to the Appendix.

\section{Statistical divergences}
\noindent We first set some general definitions and notations. Let
$P$ be some p.m. on the measurable space
$(\mathbb{R}^m,\mathcal{B}(\mathbb{R}^m))$.  Denote by $M$ the
space of all signed finite measures (s.f.m.) on
$(\mathbb{R}^m,\mathcal{B}(\mathbb{R}^m))$. Let $\varphi$ be a
convex function from $\mathbb{R}$ onto $[0,+\infty]$ with
$\varphi(1)=0$, and such that its domain, $\text{dom}\varphi
:=\left\{x\in\mathbb{R} \text{ such that }
\varphi(x)<\infty\right\}=:(a,b),$ is an interval, with endpoints
$a<1<b$, which may be bounded or unbounded, open or not. We assume
that $\varphi$ is closed\footnote{The closedness of $\varphi$
means that if $a$ or $b$ are finite then $\varphi(x)\to\varphi(a)$
when $x \downarrow a$, and $\varphi(x)\to\varphi(b)$
 when $x \uparrow b$. Note that, this is equivalent to the fact that the level sets
  $\{x\in\mathbb{R};~ \varphi(x)\leq \alpha \}$, $\forall \alpha\in\mathbb{R}$, are closed
  in $\mathbb{R}$ endowed with the usual topology.}. For any
s.f.m. $Q\in M$, the $\varphi$-divergence between $Q$ and the p.m.
$P$, when $Q$ is absolutely continuous with respect to (a.c.w.r.t)
$P$, is defined through
\begin{equation}
D_{\varphi}(Q,P):=\int_{\mathbb{R}^m}\varphi\left(\frac{dQ}{dP}(x)\right)
~dP(x),\label{divRusch}
\end{equation}
in which $\frac{dQ}{dP}(\cdot)$ denotes the Radon-Nikodym
derivative. When $Q$ is not a.c.w.r.t. $P$, we set
$D_{\varphi}(Q,P):=+\infty$. For any p.m. $P$, the mapping $Q\in
M\mapsto D_{\varphi}(Q,P)$ is convex and takes nonnegative values.
When $Q=P$ then $D_{\varphi}(Q,P)=0$. Furthermore, if the function
$x\mapsto\varphi(x)$ is strictly convex on a neighborhood of
$x=1$, then
\begin{equation}
D_{\varphi}(Q,P)=0~\text{ if and only if }~Q=P.\label{p.f.}
\end{equation}
\noindent All the above properties are presented in
\cite{Csiszar1963}, \cite{Csiszar1967} and in Chapter 1 of
\cite{Liese-Vajda1987}, for $\varphi -$divergences defined on the
set of all p.m.'s $M^{1}$. When the $\varphi $-divergences are
extended to $M$, then the same arguments as developed on $M^{1}$
hold. When defined on $M^{1}$, the Kullback-Leibler $(KL)$,
modified Kullback-Leibler $(KL_{m})$, $\chi^{2}$, modified
$\chi^{2}$ $(\chi_{m}^{2})$, Hellinger $(H)$, and $L^{1}$
divergences are respectively associated to the convex functions
$\varphi(x)=x\log x-x+1$, $\varphi(x)=-\log x+x-1$, $\varphi
(x)=\frac{1}{2}{(x-1)}^{2}$, $\varphi(x)=\frac{1}{2}{(x-1)}^{2}/x
$, $\varphi(x)=2{(\sqrt{x}-1)}^{2}$ and $\varphi(x)=\left\vert
x-1\right\vert $. All these divergences except the $L^{1}$ one,
belong to the class of the so called power divergences introduced
in \cite{Cressie-Read1984} (see also \cite{Liese-Vajda1987} and
\cite{PardoLeandro2006}). They are defined through the class of
convex functions
\begin{equation}
x\in\mathbb{R}_{+}^{\ast}\mapsto\varphi_{\gamma}(x):=\frac{x^{\gamma}-\gamma
x+\gamma-1}{\gamma(\gamma-1)}\label{gamma convex functions}
\end{equation}
if $\gamma\in\mathbb{R}\setminus\left\{  0,1\right\}  $, $\varphi
_{0}(x):=-\log x+x-1$ and $\varphi_{1}(x):=x\log x-x+1$. So, the
$KL-$divergence is associated to $\varphi_{1}$, the $KL_{m}$ to
$\varphi_{0}$, the $\chi^{2}$ to $\varphi_{2}$, the $\chi_{m}^{2}$
to $\varphi_{-1}$ and the Hellinger distance to $\varphi_{1/2}$.
We extend the definition of the power divergences functions $Q\in
M^{1}\mapsto D_{\varphi_{\gamma}}(Q,P)$ onto the whole set of
signed finite measures $M$ as follows. When the function
$x\mapsto\varphi_{\gamma}(x)$ is not defined on $]-\infty,0[$ or
when $\varphi_{\gamma}$ is defined on $\mathbb{R}$ but is not
convex,  we extend the definition of $\varphi_{\gamma}$ as follows
\begin{equation}
x\in\mathbb{R}\mapsto\varphi_{\gamma}(x)\mathds{1}_{[0,+\infty[}(x)+(+\infty
)\mathds{1}_{]-\infty,0[}(x).\label{gamma convex functions sur R}
\end{equation}
Note that for $\chi^2$-divergence, the corresponding $\varphi$
function $\varphi(x)=\frac{1}{2}(x-1)^2$ is convex and defined on
whole $\mathbb{R}$. In this paper, for technical considerations,
we assume that the  functions  $\varphi$ are strictly convex on
their domain $(a,b)$, twice continuously differentiable on
$]a,b[$, the interior of their domain. Hence, $\varphi'(1)=0$, and
for all $x\in ]a,b[$, $\varphi''(x)>0$. Here, $\varphi'$ and
$\varphi''$ are used to denote respectively the first and the
second derivative functions  of $\varphi$. Moreover, we assume
that $\varphi$ is ``essentially smooth'' in the sense that
$\lim_{x\downarrow a}\varphi'(x) =-\infty$ if $a$ is finite and
$\lim_{x\uparrow b}\varphi'(x) =+\infty$ if $b$ is finite.  Note
that the above assumptions on $\varphi$ are not restrictive, and
that all the power functions $\varphi_{\gamma}$, see (\ref{gamma
convex functions sur R}), satisfy the above conditions, including
all standard divergences.

\begin{definition}
Let $\Omega$ be some subset of $M$. The $\varphi-$divergence
between the set $\Omega$ and a p.m. $P$ is defined by
\[
D_{\varphi}(\Omega,P):=\inf_{Q\in\Omega}D_{\varphi}(Q,P).
\]
 A finite measure
$Q^{\ast}\in\Omega$, such that $D_{\varphi}(Q^{\ast},P)<\infty$
and
\[
D_{\varphi}(Q^{\ast},P)\leq D_{\varphi}(Q,P)~\text{ for all
}~Q\in\Omega,
\]
is called a projection of $P$ on $\Omega$. This projection may not
exist, or may be not defined uniquely.
\end{definition}

\section{Minimum divergence estimates}
 \noindent Let $X_{1},...,X_{n}$
denote an i.i.d. sample of a random vector $X\in\mathbb{R}^m$ with
distribution $P_{0}$. Let $P_{n}$ be the empirical measure
pertaining to this sample, namely
\[
P_{n}:=\frac{1}{n}\sum_{i=1}^{n}\delta_{X_{i}},
\]
where $\delta_{x}$ denotes the Dirac measure at point $x$, for all
$x$. We will endow our statistical approach in the global context
of s.f.m's with total mass $1$ satisfying $l$ linear constraints:
\begin{equation}
\mathcal{M}_{\theta}:=\left\{  Q\in M~\text{ such that
}~\int_{\mathbb{R}^m} dQ(x)=1\text{ and }\int_{\mathbb{R}^m}
g(x,\theta)~dQ(x)=0\right\} \label{modele signees finies}
\end{equation}
and
\begin{equation}\label{modele msf}
\mathcal{M}:=\bigcup_{\theta\in\Theta}\mathcal{M}_{\theta},
\end{equation}
sets of signed finite measures that replace
$\mathcal{M}_{\theta}^{1}$ and $\mathcal{M}^{1}$. Enhancing the
model (\ref{modele}) to the above one (\ref{modele msf}) bears a
number of improvements upon existing results; this is argued at
the end of the present Section; see also Remark \ref{remark
calcule des estim} below. The ``plug-in'' estimate of
$D_{\varphi}(\mathcal{M}_{\theta},P_{0})$ is
\begin{equation}
\widehat{D}_{\varphi}(\mathcal{M}_{\theta},P_{0}):=\inf_{Q\in\mathcal{M}_{\theta
}}D_{\varphi}(Q,P_{n})=\inf_{Q\in\mathcal{M}_{\theta}}\int_{\mathbb{R}^m}\varphi\left(
\frac {dQ}{dP_{n}}(x)\right)~dP_{n}(x).\label{plug-in est of phi}
\end{equation}
If the projection $Q^{(n)}_\theta$ of $P_n$ on
$\mathcal{M}_\theta$ exists, then it is clear that
$Q^{(n)}_\theta$ is a s.f.m. (or possibly a p.m.) a.c.w.r.t.
$P_n$; this means that the support of $Q^{(n)}_\theta$ must be
included in the set $\left\{X_1,\ldots,X_n\right\}$. So, define
the sets
\begin{equation}\label{m theta n}
\mathcal{M}_{\theta}^{(n)}:=\left\{ Q\in M~|~Q \text{ a.c.w.r.t. }
P_{n},~\sum_{i=1}^{n}Q(X_{i})=1 \text{ and } \sum_{i=1}^{n}Q(X_{i}
)g(X_{i},\theta)=0\right\},
\end{equation}
which may be seen as subsets of $\mathbb{R}^{n}$. Then, the
plug-in estimate (\ref{plug-in est of phi}) can be written as
\begin{equation}
\widehat{D}_{\varphi}(\mathcal{M}_{\theta},P_{0})=\inf_{Q\in\mathcal{M}_{\theta
}^{(n)}}\frac{1}{n}\sum_{i=1}^{n}\varphi\left(
nQ(X_{i})\right).\label{estim de phiMtheta cas cont}
\end{equation}
In the same way,
$D_{\varphi}(\mathcal{M},P_{0}):=\inf_{\theta\in\Theta}
\inf_{Q\in\mathcal{M}_{\theta}}D_{\varphi}(Q,P_{0})$ can be
estimated by
\begin{equation}
\widehat{D}_{\varphi}(\mathcal{M},P_{0}):=\inf_{\theta\in\Theta}\inf
_{Q\in\mathcal{M}_{\theta}^{(n)}}\frac{1}{n}\sum_{i=1}^{n}\varphi\left(
nQ(X_{i})\right).\label{estim de phiM  cas cont}
\end{equation}
By uniqueness of
$\arg\inf_{\theta\in\Theta}D_{\varphi}(\mathcal{M}_{\theta},P_{0})$
and since the infimum is reached at $\theta=\theta_{0}$ under the
model, we estimate $\theta_{0}$ through
\begin{equation}\label{estim de theta0 cas cont}
\widehat{\theta}_{\varphi}
:=\arg\inf_{\theta\in\Theta}\inf_{Q\in\mathcal{M}_{\theta}^{(n)}}
\frac{1}{n}\sum_{i=1}^{n}\varphi\left( nQ(X_{i})\right).
\end{equation}

\noindent Enhancing $\mathcal{M}^{1}$ to $\mathcal{M}$ and
accordingly extensions in the definitions of the $\varphi$
functions on $\mathbb{R}$ and  the $\varphi $-divergences on the
whole space of s.f.m's $M$, is motivated by the following
arguments:
\begin{enumerate}
\item [-] If the domain $(a,b)$ of the function $\varphi$ is
included in $[0,+\infty[$ then minimizing over $\mathcal{M}^{1}$
or over $\mathcal{M}$ leads to the same estimates and test
statistics. It is the case of the  $KL_m$, $KL$, modified $\chi^2$
and Hellinger divergences.

\item[-] Let $\theta$ be a given value in $\Theta$. Denote
$Q^{(1,n)}_{\theta}$ and $Q^{(n)}_\theta$, respectively, the
projection of $P_{n}$ on $\mathcal{M}_{\theta}^{1}$ and on
$\mathcal{M}_{\theta}$. If $Q^{(1,n)}_{\theta}$ satisfies
$0<Q^{(1,n)}_\theta(X_i)<1$, for all $i=1,\ldots,n,$ then
$Q^{(1,n)}_{\theta}=Q^{(n)}_\theta$. Therefore, in this case, both
approaches leads also to the same estimates and test statistics.

\item[-] It may occur that for some $\theta$ in $\Theta$ and some
$i=1,\ldots,n,$  $Q^{(1,n)}_\theta(X_i)$ is a boundary value of
$[0,1]$, hence the first order conditions are not met which makes
a real difficulty for the calculation of the estimates over the
sets of p.m. $\mathcal{M}^1_\theta$ and $\mathcal{M}^1$. However,
when $\mathcal{M}^{1}$ is replaced by $\mathcal{M}$, then this
problem does not hold any longer in particular when
$\text{dom}\varphi=\mathbb{R}$, which is the case for the
$\chi^2$-divergence. Other arguments are given in Remark
\ref{remark calcule des estim} below.
\end{enumerate}

\noindent The empirical likelihood paradigm (see \cite{Owen1988},
\cite{Owen1990}, \cite{Qin-Lawless1994} and \cite{Owen2001}),
 enters as a special case of the statistical issues related
 to estimation and tests based on
$\varphi-$divergences with $\varphi(x)=\varphi_{0}(x):=-\log
x+x-1$, namely on $KL_{m}-$divergence. Indeed, it is
straightforward to see that the empirical log-likelihood ratio
statistic for testing $\mathcal{H}_0 : P_0\in\mathcal{M}$ against
$\mathcal{H}_1 : P_0\notin\mathcal{M}$, in the context of
$\varphi$-divergences, can be written as
$2n\widehat{D}_{KL_m}(\mathcal{M},P_{0})$; and that the EL
estimate of $\theta_0$ can be written as
$\widehat{\theta}_{KL_m}=\arg\inf_{\theta\in
\Theta}\widehat{D}_{KL_m}(\mathcal{M}_\theta,P_0)$; see Remark
\ref{remarque c_0} below. In the case of the power functions
$\varphi=\varphi_{\gamma}$, the corresponding estimates
(\ref{estim de theta0 cas cont}) belong to the class of GEL
estimates introduced by \cite{Smith1997} and
\cite{NeweySmith2004}, and (\ref{estim de phiMtheta cas cont}) in
 this case  are the empirical Cressie-Read statistics
introduced by \cite{Baggerly1998}
and \cite{Corcoran1998}; see Remark \ref{remark on GEL estimates} below.\\

\noindent The constrained optimization problems (\ref{estim de
phiMtheta cas cont}), (\ref{estim de phiM cas cont}) and
(\ref{estim de theta0 cas cont}) can be transformed into
unconstrained ones making use of some arguments of ``duality''
which we briefly state below from \cite{Rockafellar1970}.  On the
other hand, the obtaining of asymptotic statistical results of the
estimates and the test statistics, under misspecification or under
alternative hypotheses, requires  handle existence conditions and
characterization of the projection of $P_{0}$ on the submodel
$\mathcal{M}_{\theta}$ or on the model $\mathcal{M}.$ This also
will be considered through duality, along the following Section.

\section{Dual representation of $\varphi-$divergences under
constraints} \noindent This Section is central for our purposes.
Indeed, it provides the explicit form of the proposed estimates by
transforming the constrained problems (\ref{estim de phiMtheta cas
cont}) to unconstrained ones, using Lagrangian duality  which is a
classical tool in optimization theory. This Section adapts this
formalism to the context of divergences and the present
statistical setting. The Lagrangian ``dual'' problem,
corresponding to the ``primal'' one
\begin{equation}
\inf_{Q\in\mathcal{M}_{\theta}}D_\varphi(Q,P_{0})\label{Primal
theo}
\end{equation}
and its empirical counterpart (\ref{estim de phiMtheta cas cont}),
make use of the so-called Fenchel-Legendre transform  of
$\varphi$, defined through
\begin{equation}
\psi : t\in \mathbb{R}\mapsto
\psi(t):=\sup_{x\in\mathbb{R}}\left\{tx-\varphi(x)\right\}.
\end{equation}
The ``dual'' problems associated to $(\ref{Primal theo})$  and
(\ref{estim de phiMtheta cas cont}) are respectively
\begin{equation}\label{dual theo}
\sup_{t\in\mathbb{R}^{1+l}}\left\{  t_{0}-\int_{\mathbb{R}^m} \psi
( t_{0}+\sum_{j=1}^{l}t_{j}g_{j}(x,\theta))~ dP_{0}(x)\right\},
\end{equation}
and
 \begin{equation}\label{dual empi}
\sup_{t\in\mathbb{R}^{1+l}}\left\{
t_{0}-\frac{1}{n}\sum_{i=1}^{n}\psi(
t_{0}+\sum_{j=1}^{l}t_{j}g_{j}(X_{i},\theta))\right\}.
\end{equation}
In the following Propositions \ref{proposition 1} and
\ref{proposition 2}, we state sufficient conditions under which
the primal problems (\ref{Primal theo}) and (\ref{estim de
phiMtheta cas cont}) coincide respectively with the dual ones
(\ref{dual theo}) and (\ref{dual empi}). First, recall some
properties of the convex conjugate $\psi$ of $\varphi$. For the
proofs, we can refer to Section 26 in \cite{Rockafellar1970}. The
function $\psi$ is convex and closed, its domain is an interval
with endpoints
\begin{equation}\label{domain de varphi*}
a^*:=\lim_{x\to  -\infty}\frac{\varphi(x)}{x},\quad b^*:=\lim_{x
\to +\infty}\frac{\varphi(x)}{x}
\end{equation}
satisfying $a^*<0<b^*$ with $\psi(0)=0$.  The strict convexity of
$\varphi$ on its domain $(a,b)$ is equivalent to the condition
that its conjugate $\psi$ is essentially smooth, i.e.,
differentiable with
\begin{equation}
\begin{array}
[c]{ccccc} \lim_{t\downarrow a^*}{\psi}^{\prime}(t) & = & -\infty
&
\text{ if } & a^* \text{ is finite},\\
\lim_{t\uparrow b^*}{\psi}^{\prime}(t) & = & +\infty & \text{ if }
& b^* \text{ is finite}.
\end{array}
\end{equation}
Conversely, $\varphi$ is essentially smooth on its domain $(a,b)$
if and only if $\psi$ is strictly convex on its domain
$(a^*,b^*)$. In all the sequel, we assume additionally that
$\varphi$ is essentially smooth. Hence, $\psi$ is strictly convex
on its domain $(a^*,b^*)$, and it holds that
 $$a^*=\lim_{x\downarrow a}\varphi'(x),\quad \quad b^*=
 \lim_{x\uparrow b}\varphi'(x),$$
 and
\begin{equation}\label{forme explicite de psi}
\psi(t)=t{\varphi'}^{-1}(t)-\varphi\left(
{\varphi'}^{-1}(t)\right),\quad \text{for all } t\in ]a^*,b^*[,
\end{equation}
where ${\varphi'}^{-1}$ denotes the inverse function of
$\varphi'$. It holds also that $\psi$ is twice continuously
differentiable on $]a^*,b^*[$ with
\begin{equation}\label{derivee de psi}
 \psi'(t)={\varphi'}^{-1}(t) \quad \text{and}\quad
 \psi''(t)=\frac{1}{\varphi''\left({\varphi'}^{-1}(t)\right)}.
\end{equation}
In particular, $\psi'(0)=1$ and $\psi''(0)=1$.  Obviously, since
$\varphi$ is assumed to be closed, we have
 $$\varphi(a)=\lim_{x\downarrow a}\varphi(x)\quad \text{ and }\quad
 \varphi(b)=\lim_{x\uparrow b}\varphi(x),$$
 which may be finite or infinite. Hence, by closedness of $\psi$,
 we have
 $$\psi(a^*)=\lim_{t\downarrow a^*}\psi(x)\quad \text{ and }\quad
 \psi(b^*)=\lim_{t\uparrow b^*}\psi(t).$$
 Finally, the first and second derivatives of $\varphi$ in $a$ and $b$  are
 defined to be the limits of $\varphi'(x)$ and $\varphi''(x)$ when $x\downarrow a$ and
 when $x\uparrow b$. The first and second derivatives of $\psi$ in $a^*$ and $b^*$
 are defined in
 a similar way. In Table \ref{table convex conjugates}, we give
 the convex conjugates $\psi$ of some  standard functions $\varphi$, associated to
 some standard divergences. We determine also their domains, $(a,b)$ and
 $(a^*,b^*)$.

\begin{table}[h]
\caption{Convex conjugates for some standard divergences.}
\begin{tabular}{|l||l|l||l|l|}
  \hline
  $D_\varphi$ & $\varphi$ & $\text{dom}\varphi $ & $\text{dom}\psi $ & $\psi$ \\
  \hline
  \hline
  $D_{KL_m}$ &  $\varphi(x):=-\log x +x -1$ & $]0,+\infty[$  & $]-\infty,1[$ & $\psi(t)= - \log(1-t)$ \\
  \hline
  $D_{KL}$ &  $\varphi(x):=x\log x -x +1$ & $[0,+\infty[$ & $\mathbb{R}$ & $\psi(t)=e^t-1$  \\
  \hline
  $D_{\chi^2_m}$ &  $\varphi(x):= \frac{1}{2}\frac{\left(x-1\right)^2}{x}$
  & $]0,+\infty[$ &
  $\left]-\infty,\frac{1}{2}\right]$ & $\psi(t)=1-\sqrt{1-2t}$ \\
  \hline
  $D_{\chi^2}$ &  $\varphi(x):= \frac{1}{2}\left(x-1\right)^2$ &  $\mathbb{R}$ &  $\mathbb{R}$ & $\psi(t)=\frac{1}{2}t^2+t$ \\
  \hline
  $D_H$ &  $\varphi(x):=2(\sqrt{x}-1)^2$ & $[0,+\infty[$ &  $]-\infty, 2[$ & $\psi(t)=\frac{2t}{2-t}$ \\
  \hline
  $D_{\varphi_\gamma}$ &  $\varphi(x):=\frac{x^\gamma -\gamma x +\gamma -1}{\gamma (\gamma -1)}$ & $--$  & $--$ &
  $\psi(t)=\frac{1}{\gamma}\left(\gamma t -t+1\right)^{\frac{\gamma}{\gamma -1}}-\frac{1}{\gamma}$
  \\
  \hline
\end{tabular}
\label{table convex conjugates}
\end{table}

\begin{proposition}\label{proposition 1}
\label{Propo carct empirique} Let $\theta$ be a given value in
$\Theta$. If there exists $Q_{0}$ in $\mathcal{M}_{\theta}^{(n)}$
such that
\begin{equation}\label{cond Owen}
a< Q_{0}(X_{i})<b, \quad \text{for all} \quad i=1,\ldots,n,
\end{equation}
then
\begin{equation}\label{egalite duale 1}
\inf_{Q\in\mathcal{M}_{\theta}^{(n)}}D_{\varphi}(Q,P_{n})=\sup_{t\in
\mathbb{R}^{1+l}}\left\{ t_{0}-\frac{1}{n}\sum_{i=1}^{n}\psi(
t_{0} +\sum_{j=1}^{l}t_{j}g_{j}(X_{i},\theta)) \right\}
\end{equation}
with dual attainment. Conversely,  if there exists some dual
optimal solution
$\widehat{t}:=(\widehat{t}_0,\widehat{t}_1,\ldots,\widehat{t}_l)^{\top}\in\mathbb{R}^{1+l}$
such that
\begin{equation}\label{cond Owen 2}
 a^*< \widehat{t_0}+\sum_{j=1}^l \widehat{t_j}g_j(X_i,\theta)<b^*, \quad \text{for
 all}\quad
 i=1,\ldots,n,
\end{equation}
then the equality (\ref{egalite duale 1}) holds, and  the unique
optimal solution of the primal problem
$\inf_{Q\in\mathcal{M}_{\theta}^{(n)}}D_{\varphi}(Q,P_{n})$,
namely the projection of $P_n$ on $\mathcal{M}^{(n)}_\theta$, is
given by
\[
Q^{(n)}_\theta(X_i)=\frac{1}{n}
{\varphi'}^{-1}(\widehat{t_0}+\sum_{j=1}^{l}\widehat{t_j} g_j(
X_i,\theta)),\quad i=1,\ldots,n,
\]
where
$\widehat{t}:=(\widehat{t}_0,\widehat{t}_1,\ldots,\widehat{t}_l)^{\top}$
is solution of the system of equations
\begin{equation*}
\left\{\begin{array}{lll}
  1-\frac{1}{n}\sum_{i=1}^{n} {\varphi'}^{-1}(\widehat{t_0}
+\sum_{j=1}^l\widehat{t_j} g_j (X_i,\theta)) & = & 0, \\
  -\frac{1}{n}\sum_{i=1}^n g_{j}(X_i,\theta)
 {\varphi'}^{-1}( \widehat{t_0}+\sum_{j=1}^{l}\widehat{t_j}g_{j}(
X_{i},\theta)) & = & 0, \quad j=1,\ldots,l.\\
\end{array} \right.
\end{equation*}\\
\end{proposition}

\begin{remark}{\rm
For the $\chi^{2}-$divergence, we have $a=-\infty$ and
$b=+\infty$. Hence, condition (\ref{cond Owen}) holds whenever
$\mathcal{M}_{\theta}^{(n)}$ is not void. More generally, the
above Proposition holds for any $\varphi$-divergence with
 $\text{dom}\varphi=\mathbb{R}$.}
\end{remark}

\begin{remark}{\rm
Assume that $g(x,\theta):=(x-\theta)^\top$. So, for any divergence
$D_\varphi$ with $\text{dom}\varphi = ]0,+\infty[$, which is the
case of the modified $\chi^2$ divergence and the modified
Kullback-Leibler divergence (or equivalently EL method), condition
$(\ref{cond Owen})$ means that $\theta$ is an interior point of
the convex hull of the data $(X_{1},...,X_{n})$. This is precisely
what is checked in \cite{Owen1990}, p. 100, for the EL method; see
also \cite{Owen2001}.}
\end{remark}

\noindent For the asymptotic counterpart of the above results we
have; see Theorem 1 in \cite{Bronia_Kez2006_STUDIA}:
\begin{proposition}\label{proposition 2}
\label{Prop caract theo} Let $\theta$ be a given value in
$\Theta$. Assume that $\int |g_j(x,\theta)|~dP_0(x)<\infty$, for
all $j=1,\ldots,l$. If there exists $Q_{0}$ in
$\mathcal{M}_{\theta}$ with $D_\varphi(Q_0,P_0)<\infty$
and\footnote{The strict inequalities  (\ref{condition de
qualification}) mean that
$P_0\left\{x\in\mathbb{R}^m~|~\frac{dQ_0}{dP_0}(x)\leq
a\right\}=P_0\left\{x\in\mathbb{R}^m~|~\frac{dQ_0}{dP_0}(x)\geq
b\right\}=0.$}
\begin{equation}\label{condition de qualification}
a<\inf_{x\in\mathbb{R}^m}\frac{dQ_{0}}{dP_{0}}(x)\leq\sup_{x\in\mathbb{R}^m}\frac{dQ_{0}}{dP_{0}
}(x)<b, \quad P_{0}-a.s.,
\end{equation}
then
\begin{equation}\label{egalite duale 2}
\inf_{Q\in\mathcal{M}_{\theta}}D_{\varphi}(Q,P_{0})=\sup_{t\in\mathbb{R}^{1+l}
}\left\{  t_{0}-\int_{\mathbb{R}^m} \psi(
t_{0}+\sum_{j=1}^{l}t_{j}g_{j}( x,\theta))~dP_{0}(x)\right\}
\end{equation}
 with dual attainment. Conversely, if  there exists some dual optimal solution
$t^*$ which is an interior point of the set
\begin{equation}\label{condition de qualification 2}
\left\{t\in\mathbb{R}^{1+l}~\text{such that}~ \int_{\mathbb{R}^m}
|\psi(t_{0}+\sum_{j=1}^{l}t_jg_{j}(x,\theta))| ~dP_0(x) <\infty
\right\},
\end{equation}
 then the dual equality (\ref{egalite duale 2}) holds, and
 the unique optimal solution $Q^*_\theta$  of the
primal problem
$\inf_{Q\in\mathcal{M}_{\theta}}D_{\varphi}(Q,P_{0})$, namely the
projection of $P_0$ on $\mathcal{M}_\theta$, is given by
\begin{equation*}
\frac{dQ^*_\theta}{dP_{0}}(x)={\varphi'}^{-1}(t_{0}^*+\sum_{j=1}^{l}t_j^*g_{j}(x,\theta)),
\end{equation*}
where $t^*:=(t^*_0,t^*_1,\ldots,t^*_l)^{\top}$ is solution of the
system of equations
\begin{equation}\label{systeme}
\left\{\begin{array}{lll}
  1-\int {\varphi'}^{-1}( t_{0}^*+\sum_{j=1}^{l}t_j^*g_{j}( x,\theta)
)~  dP_{0}(x) & = & 0, \\
  -\int g_{j}(x,\theta)
{\varphi'}^{-1}(t_{0}^*+\sum_{j=1}^{l}t_j^*g_{j}(x,\theta))~
dP_{0}(x) & = & 0,\quad j=1,\ldots,l. \\
\end{array} \right.
\end{equation}
Furthermore, $t^*$ is unique if the functions
$\mathds{1}_{\mathbb{R}^m}, g_1(.,\theta),\ldots,g_l(.,\theta)$
are linearly independent in the sense that
$P_0\left\{x\in\mathbb{R}^m~|~t_0+\sum_{j=1}^lt_jg_j(x,\theta)\neq
0\right\}>0$ for all
$t\in\mathbb{R}^m$ with $t\neq 0.$\\
\end{proposition}
\noindent For sake of brevity and clearness, we must introduce
some additional notations. In all the sequel, $\|x\|$ denotes the
norm of $x$ defined by  $\|x\|:=\sup_{i}|x_i|$ for  any vector
$x:=(x_1,\ldots,x_k)^\top\in\mathbb{R}^k$, and for any matrix $A$,
the norm of $A$ is defined by $\|A\|:=\sup_{i,j}|a_{i,j}|$. Denote
by $\overline{g}$ the vector valued function
$\overline{g}:=(\mathds{1}_{\mathbb{R}^m},g_1,\ldots,g_l)^{\top}\in\mathbb{R}^{1+l}$.
For any p.m. $P$ on $(\mathbb{R}^m,\mathcal{B}(\mathbb{R}^m))$ and
any real  measurable function $f$ from
$(\mathbb{R}^m,\mathcal{B}(\mathbb{R}^m))$ to
$(\mathbb{R},\mathcal{B}(\mathbb{R}))$, denote
$$Pf:=\int_{\mathbb{R}^m}f(x)~dP(x).$$
 Let
 $$t^{\top} \overline{g}(x,\theta):=t_0+\sum_{j=1}^lt_jg_j(x,\theta)$$
 and
 \begin{equation}\label{m x theta t} m(x,\theta,t):= t_0-
 \psi(t^{\top} \overline{g}(x,\theta)), \text{
for all } x\in\mathbb{R}^m,
\theta\in\Theta\subset\mathbb{R}^d,t\in \mathbb{R}^{1+l}.
  \end{equation}
\noindent Note that the $\sup$ in (\ref{egalite duale 1}) and
(\ref{egalite duale 2}) can be restricted, respectively, to the
sets
 \begin{equation}\label{lambda n theta} \Lambda_\theta^{(n)}:=\left\{t\in\mathbb{R}^{1+l}~|~
a^*<t^{\top}\overline{g}(X_i,\theta)<b^*, ~\text{ for all }
i=1,\ldots,n\right\}\end{equation} and
 \begin{equation}\label{lambda theta}\Lambda_\theta:=\left\{t\in\mathbb{R}^{1+l}~|~
 \int_{\mathbb{R}^m} |\psi(t_{0}+\sum_{j=1}^{l}t_jg_{j}(x,\theta))|
~dP_0(x) <\infty \right\}.\end{equation}
 In view of
the above two Propositions \ref{proposition 1} and
\ref{proposition 2}, we redefine the estimates (\ref{estim de
phiMtheta cas cont}), (\ref{estim de phiM cas cont}) and
(\ref{estim de theta0 cas cont}) as follows
\begin{equation}\label{estimateur1}
\widehat{D}_\varphi\left(\mathcal{M}_\theta,P_0\right):=
\sup_{t\in\Lambda_\theta^{(n)}}\frac{1}{n}\sum_{i=1}^n
m(X_i,\theta,t):=\sup_{t\in\Lambda_\theta^{(n)}} P_nm(\theta,t),
\end{equation}
\begin{equation}\label{estimateur2}
\widehat{D}_\varphi\left(\mathcal{M},P_0\right):=
\inf_{\theta\in\Theta}\sup_{t\in\Lambda_\theta^{(n)}}\frac{1}{n}\sum_{i=1}^n
m(X_i,\theta,t):=\inf_{\theta\in\Theta}\sup_{t\in\Lambda_\theta^{(n)}}
P_nm(\theta,t)
\end{equation}
and
\begin{equation}\label{estimateur3}
\widehat{\theta}_\varphi:=
\arg\inf_{\theta\in\Theta}\sup_{t\in\Lambda_\theta^{(n)}}\frac{1}{n}\sum_{i=1}^n
m(X_i,\theta,t):=\arg\inf_{\theta\in\Theta}\sup_{t\in\Lambda_\theta^{(n)}}P_nm(\theta,t).
\end{equation}

\begin{remark}\label{remarque c_0}{\rm
 When $\varphi(x)=-\log x+x-1$, then the estimate (\ref{estim de theta0 cas cont})
 clearly coincides with the EL one, so
 it can be seen as the value of the parameter which
 minimizes the $KL_m$-divergence between the model $\mathcal{M}$ and the empirical measure
 $P_n$ of the data $X_1,\ldots, X_n$. The statistic $2n\widehat{D}_{KL_m}(\mathcal{M},P_0)$, see (\ref{estim de phiM  cas
 cont}), coincides with the empirical likelihood ratio statistic  associated
 to the null hypothesis $\mathcal{H}_0:P_0\in\mathcal{M}$ against the alternative $\mathcal{H}_1:P_0\not\in\mathcal{M}$.
  The  dual representation of $\widehat{D}_{KL_m}(\mathcal{M},P_0)$, see (\ref{estimateur2}) and (\ref{egalite duale 1}),  is
  $$\widehat{D}_{KL_m}(\mathcal{M},P_0)=\inf_{\theta\in\Theta}\sup_{t\in\Lambda_\theta^{(n)}}\left\{t_0+\frac{1}{n}\sum_{i=1}^n
  \log(1-t_0-
  \sum_{j=1}^l t_jg_j(X_i,\theta))\right\}.$$
  For $a$ given $\theta\in\Theta$, the $KL_m$-projection $Q^{(n)}_\theta$, of $P_n$
  on $\mathcal{M}_{\theta}$, is given by (see Proposition \ref{proposition 1})
   $$\frac{1}{Q^{(n)}_\theta(X_i)}=n\left(1-\widehat{t}_0-\sum_{j=1}^l\widehat{t}_jg(X_i,\theta)\right),\quad i=1,\dots,n,$$
  which, multiplying by $Q^{(n)}_\theta(X_i)$ and summing upon $i=1,\ldots,n,$ yields
  $\widehat{t}_0=0$. Therefore, $t_0$ can be omitted, and the above representation can be rewritten
  as follows
  $$\widehat{D}_{KL_m}(\mathcal{M},P_0)=\inf_{\theta\in\Theta}\sup_{t_1,\ldots,t_l}\left\{\frac{1}{n}\sum_{i=1}^n
  \log(1+
  \sum_{j=1}^l t_jg_j(X_i,\theta))\right\}$$
  and then
   \begin{equation}\label{representation duale de EL}
  \widehat{\theta}_{KL_m}=\widehat{\theta}_{EL}=\arg\inf_{\theta\in\Theta}\sup_{t_1,\ldots,t_l}\left\{\frac{1}{n}\sum_{i=1}^n
  \log(1+
  \sum_{j=1}^l t_jg_j(X_i,\theta))\right\}
   \end{equation}
  in which the $\sup$ is taken over  the set $$\left\{(t_1,\ldots,t_l)^\top\in\mathbb{R}^m ~ |~ -1<
  \sum_{j=1}^lt_jg_j(X_i,\theta)<+\infty, \text{ for all }
  i=1,\ldots,n
  \right\}.$$
  The formula (\ref{representation duale de EL}) is the ordinary dual representation of the EL estimate;
  see \cite{Qin-Lawless1994} and \cite{Owen2001}.}\\
\end{remark}

\begin{remark}\label{remark on GEL estimates}{\rm Consider the power divergences, associated to the power
functions $\varphi_\gamma$; see (\ref{gamma convex functions}) and
(\ref{gamma convex functions sur R}). We will show that the
estimates $\widehat{\theta}_{\varphi_{\gamma}}$ belong to the
class of GEL estimators introduced by \cite{Smith1997} and
\cite{NeweySmith2004}. The projection $Q^{(n)}_\theta$ of $P_n$ on
$\mathcal{M}_\theta$ is given by
$$Q^{(n)}_\theta(X_i)=\left((\gamma-1) (\widehat{t}_0+\sum_{j=1}^l \widehat{t}_jg(X_i,\theta))
+1\right)^{1/(\gamma-1)},\quad i=1,\ldots,n.$$ Using the
constraint $\sum_{i=1}^n Q^{(n)}_\theta(X_i)=1$, we can explicit
$\widehat{t}_0$ in terms of $\widehat{t}_1,\ldots,\widehat{t}_l$,
and hence the $\sup$ in the dual representation
(\ref{estimateur3}) can be reduced to a subset of $\mathbb{R}^l,$
as in \cite{NeweySmith2004}. When $\varphi(x)=\frac{1}{2}(x-1)^2$,
it is straightforward to see that the corresponding estimate
$\widehat{\theta}_\varphi$ coincides with the continuous updating
estimator of
 \cite{Hansen_Healton_Yaron1996}.}
\end{remark}

\begin{remark}{\rm (\textbf{Numerical calculation of the estimates and the specific role of the
$\chi^2$-divergence}). The computation of $\widehat{t}(\theta)$
for fixed $\theta\in\Theta$ as defined in (\ref{systeme}) is
difficult when handling a generic divergence. In the particular
case of $\chi^2$-divergence, i.e., when
$\varphi(x)=\frac{1}{2}(x-1)^2$, optimizing on all s.f.m's, the
system (\ref{systeme}) is linear; we thus easily obtain an
explicit form for $\widehat{t}(\theta)$, which in turn allows for
a single gradient descent when optimizing upon $\Theta$. This
procedure is useful in order to compute the estimates for all
other divergences (for which the corresponding system is non
linear) including EL, since it provides an easy starting point for
the resulting double gradient descent. Moreover,
\cite{Hjort_McKeague_VanKeilegom2009} extend the EL approach, to
more general moment condition models, allowing the number of
constraints to increase with growing sample size. In this case,
the computation of EL estimate is more complex, and the same idea
as above can help to solve the problem. \label{remark calcule des
estim}}
\end{remark}

\section{Asymptotic properties of the estimates of the parameter and
 the divergences}

\subsection{Asymptotic properties under the model}
This Section addresses Problems 1 and 2,  aiming at testing the
null hypothesis $\mathcal{H}_0 : P_{0} \in \mathcal{M}$ against
the alternative  $\mathcal{H}_1: P_{0}\not\in \mathcal{M}$. We
derive the limiting distributions of the proposed  test statistics
which are the estimated divergences between the model
$\mathcal{M}$ and $P_0$. We also derive the limiting distributions
of the estimates of $\theta_0$. The following two Theorems
\ref{theoreme 1} and \ref{theoreme 2} extend Theorems 3.1 and 3.2
in \cite{NeweySmith2004} to the context of divergence based
approach. The Assumptions which we will consider match those of
Theorems 3.1 and 3.2 in  \cite{NeweySmith2004}.

\begin{assumption}
\label{assumption 1} a) $P_{0}\in \mathcal{M}$ and $\theta _{0}\in
\Theta $ is the unique solution of $\mathbb{E}\left[ g(X,\theta
)\right] =0$; b) $\Theta \subset \mathbb{R}^{d}$ is compact; c)
$g(X,\theta )$ is continuous at each $\theta \in \Theta $ with
probability one; d) $\mathbb{E}\left[ \sup_{\theta \in \Theta
}\Vert g(X,\theta )\Vert ^{\alpha }\right] <\infty $ for some
$\alpha >2$; e) the matrix $\Omega :=\mathbb{E}\left[ g(X,\theta
_{0})g(X,\theta _{0})^{\top}\right] $ is nonsingular.
\end{assumption}

\begin{theorem}\label{theoreme 1}
 Under Assumption \ref{assumption 1},  with probability approaching one as $n\to\infty$,
 the estimate $\widehat{\theta}_\varphi$
 exists, and converges to
 $\theta_0$ in probability. $\frac{1}{n}\sum_{i=1}^n
 g(X_i,\widehat{\theta}_\varphi)=O_P(1/\sqrt{n})$,
  $\widehat{t}(\widehat{\theta}_\varphi):=\arg\sup_{t\in\Lambda^{(n)}_{\widehat{\theta}_\varphi}}
  P_nm(\widehat{\theta}_\varphi,t)$
  exists and belongs to $\text{int}(\Lambda^{(n)}_{\widehat{\theta}_\varphi})$
  with probability approaching one  as $n\to\infty$, and
  $\widehat{t}(\widehat{\theta}_\varphi)=O_P(1/\sqrt{n})$.
\end{theorem}

\noindent In order to obtain asymptotic normality, we need some
additional Assumptions. Denote by $G$ the matrix
$G:=\mathbb{E}\left[\partial g(X,\theta_0)/\partial
\theta\right]$.

\begin{assumption}
\label{assumption 2} a) $\theta _{0}\in \text{int}(\Theta )$; b)
with probability one, $g(X,\theta )$ is continuously
differentiable in a neighborhood $N_{\theta_0}$ of $\theta _{0}$,
and $\mathbb{E}\left[ \sup_{\theta \in N_{\theta_0}}\Vert
\partial g(X,\theta )/\partial \theta \Vert \right] <\infty $; c)
$\text{rank}(G)=d.$
\end{assumption}

\begin{theorem}\label{theoreme 2}
 Assume that  Assumptions \ref{assumption 1} and \ref{assumption 2} hold.
 Then,
 \begin{enumerate}
\item [1)]
$\sqrt{n}\left(\widehat{\theta}_\varphi-\theta_0\right)$
 converges in distribution
 to a centered normal random vector with covariance matrix   $$V:=\left[G \Omega^{-1}G^{\top}\right]^{-1}.$$
  \item  [2)] If $l>d$, then the statistic $2n\widehat{D}_\varphi(\mathcal{M},P_0)$
  converges in distribution to a $\chi^2$ random variable with
  $(l-d)$ degrees of freedom.
\end{enumerate}
\end{theorem}

\begin{remark} \label{remarque test du model}{\rm The above Theorem allows to perform statistical tests (of the model)
with asymptotic level $\alpha\in ]0,1[$. Consider the null
hypothesis
 \begin{equation}\label{test de model}\mathcal{H}_0 : P_0\in\mathcal{M} \quad
 \text{against the alternative}\quad \mathcal{H}_1:P_0\not\in\mathcal{M}.
 \end{equation}
 The critical region is then
 $$C_\varphi:=\left\{2n\widehat{D}_\varphi(\mathcal{M},P_0)>q_{(1-\alpha)}\right\}$$
 where $q_{(1-\alpha)}$ is the $(1-\alpha)$-quantile of the
 $\chi^2(l-d)$ distribution. When $\varphi(x)=-\log x +x -1$, it is straightforward to see that the
 corresponding test is the empirical likelihood ratio one;
 see \cite{Qin-Lawless1994}.}
\end{remark}

\subsection{Asymptotic properties of the estimates of the divergences for a
given value of the parameter}

For a given $\theta \in\Theta$, consider the test problem of the
null hypothesis  $\mathcal{H}_0: P_{0}\in \mathcal{M}_{\theta}$
against two different families of alternative hypotheses:
$\mathcal{H}_1: P_{0}\notin \mathcal{M}_{\theta}$ and
$\mathcal{H}_1':P_{0}\in \mathcal{M}\setminus\mathcal{M}_\theta.$
Those two tests address different situations since $\mathcal{H}_1$
may include misspecification of the model. We give two different
test statistics each pertaining to one of the situations and
derive their limiting distributions both under $\mathcal{H}_0$ and
under the alternatives. As a by product, we also derive confidence
areas for the true value $\theta_0$ of the parameter. We will
first state the convergence in probability of
$\widehat{D}_\varphi(\mathcal{M}_{\theta},P_0)$ to
$D_\varphi(\mathcal{M}_{\theta},P_0)$, and then we obtain the
limiting distribution of
$\widehat{D}_\varphi(\mathcal{M}_{\theta},P_0)$ both when
$P_0\in\mathcal{M}_{\theta}$ and when $P_0\not\in
\mathcal{M}_{\theta}$. Obviously, when $P_0\in
\mathcal{M}_{\theta}$, this means that $\theta=\theta_0$ since the
true value $\theta_0$ of the parameter is assumed to be unique.

\begin{assumption}
\label{assumption 3}  a) $P_0\in\mathcal{M}_{\theta}$ and $\theta$
is the unique solution of $\mathbb{E}\left[g(X,\theta)\right]=0$;
b) $\mathbb{E} \left[\|g(X,\theta)\|^\alpha\right]<\infty$ for
some $\alpha>2$; c) the matrix $\Omega:=\mathbb{E}\left[
g(X,\theta)g(X,\theta)^{\top}\right]$ is nonsingular.
\end{assumption}

\begin{theorem}\label{theoreme 3}
Under Assumption \ref{assumption 3}, we have
 \begin{enumerate}
  \item [1)]
  $\widehat{t}(\theta):=\arg\sup_{t\in\Lambda_\theta^{(n)}}P_nm(\theta,t)$
  exists and belongs to $\text{int}(\Lambda^{(n)}_\theta)$ with
  probability approaching one as $n\to\infty$, and
  $\widehat{t}(\theta)=O_P(1/\sqrt{n})$.
  \item [2)] The statistic
  $2n\widehat{D}_\varphi(\mathcal{M}_\theta,P_0)$ converges in
  distribution to a $\chi^2(l)$ random variable.
 \end{enumerate}
\end{theorem}

\noindent In order to obtain the limiting distribution of the test
statistic $2n\widehat{D}_{\varphi }\left( \mathcal{M}_{\theta
},P_{0}\right) $ under the alternative $\mathcal{H}_1:P_{0}\notin
\mathcal{M}_{\theta }$, including misspecification, the following
Assumption is needed.

\begin{assumption} \label{assumption 4}
 a) $P_{0}\not\in \mathcal{M}_{\theta}$, and
$t^{\ast}(\theta ):=\arg \sup_{t\in
\Lambda_\theta}\mathbb{E}\left[ m(X,\theta,t)\right]$ exists and
is an interior point of $\Lambda_\theta$; b) $\mathbb{E}\left[
\sup_{t\in N_{t^*(\theta)}}|m(X,\theta ,t)|\right] <\infty $ for
some compact set $N_{t^*(\theta)}\subset \Lambda_\theta$ such that
$t^{\ast}(\theta)\in \text{int}(N_{t^*(\theta)})$; c) the
functions $\mathds{1}_{\mathbb{R}^{m}},g_{1},\ldots ,g_{l}$ are
linearly independent in the sense that : $P_{0}\left\{
x\in\mathbb{R}^m~|~t_{0}+\sum_{j=1}^{l}t_{j}g_{j}(x,\theta )\neq
0\right\}
>0$ ~ for all $t\in \mathbb{R}^{1+l}$ with $t\neq 0$.

\begin{remark}\label{remark 1 sur les conditions}{\rm
 Assumption \ref{assumption 4}.c above ensures the strict concavity of the function
 $t\in \Lambda_\theta\mapsto \mathbb{E}\left[m\left( X,\theta
,t\right)\right]$ on the convex set $\Lambda_\theta$, which
implies that $t^*(\theta)$ is unique. It can be replaced by the
following Assumption : there exists a neighborhood,
$N_{t^*(\theta)}\subset\Lambda_\theta$, of $t^*(\theta)$, such
that $\mathbb{E}\left[\sup_{t\in N_{t^*(\theta)}}\left\|\partial
m(X,{\theta},t)/\partial t\right\|\right]<\infty$,
$\mathbb{E}\left[\sup_{t\in
N_{t^*(\theta)}}\left\|\partial^2m(X,{\theta},t)/\partial t
^2\right\|\right]<\infty$ and  the matrix
$\mathbb{E}\left[\partial^2m(X,{\theta},t^*(\theta))/\partial t
^2\right]$ is nonsingular; which implies also that $t^*(\theta)$
is unique.}
\end{remark}
\end{assumption}

\begin{theorem}\label{theoreme 4}
Under Assumption \ref{assumption 4}, when $P_0\not\in
\mathcal{M}_{\theta}$, we have
\begin{enumerate}
 \item [1)]  $\widehat{t}(\theta)$ converges in probability to $t^*(\theta)$.
 \item [2)]  $\widehat{D}_\varphi(\mathcal{M}_{\theta},P_0)$ converges
 in probability to $D_\varphi(\mathcal{M}_{\theta},P_0).$
\end{enumerate}
\end{theorem}

\noindent We now give the limiting distribution of the test
statistic under $\mathcal{H}_1.$ We need the following additional
condition.

\begin{assumption}
\label{assumption 5} a) There exists
$N_{t^*(\theta)}\subset\Lambda_\theta$,  some compact neighborhood
of $t^*(\theta)$, such that
$$\mathbb{E}[\sup_{t\in N_{t^*(\theta)}} \|\partial m(X,\theta,t^*(\theta))/\partial t\|]<\infty,\quad
\mathbb{E}[\sup_{t\in N_{t^*(\theta)}} \|\partial^2
m(X,\theta,t^*(\theta))/\partial t^2\|]<\infty;$$
  b) as $\delta\to 0$,
    $$\mathbb{E}\left\{\sup_{\{t; \|t-t^*(\theta)\|\leq \delta\}}\left\|\partial^2m(X,\theta,t)/\partial t^2-
    \partial^2m(X,\theta,t^*(\theta))/\partial t^2\right\|\right\}\to 0;$$
 c) $\mathbb{E}\left[m(X,{\theta},t^*(\theta))^2\right]<\infty$, $\mathbb{E}\left[\|\partial m(X,{\theta},
 t^*(\theta))/\partial t\|^2\right]<\infty$\\
 and the matrix
$\mathbb{E}\left[\partial^2m(X,{\theta},t^*(\theta))/\partial t
^2\right]$ is nonsingular.
\end{assumption}

\begin{remark}\label{remark 1 sur les conditions de 3 ordre}
 Assumption \ref{assumption 5}.b is used here to relax the
 condition on the third derivatives (in $t$) of the function $t\mapsto
 m(X,\theta,t)$.
\end{remark}

\begin{theorem}\label{theoreme 5}
 Under Assumptions \ref{assumption 4} and \ref{assumption 5}, we
 have
 \begin{enumerate}
  \item [1)] $\sqrt{n}(\widehat{t}(\theta)-t^*(\theta))$ converges in
  distribution to a centered normal random vector with covariance matrix
   $$\left[\mathbb{E}\left[m''(X,{\theta},t^*)\right] \right]^{-1}
   \mathbb{E}\left[m'(X,{\theta},t^*)m'(X,{\theta},t^*)^{\top}\right]
   \left[\mathbb{E}\left[m''(X,{\theta},t^*)\right] \right]^{-1}.$$
  \item [2)] $\sqrt{n}\left(\widehat{D}_\varphi(\mathcal{M}_{\theta},P_0)-D_\varphi(\mathcal{M}_{\theta},P_0)\right)$
  converges in distribution to a centered normal random variable
  with variance
   $$\sigma^2(\theta):=\mathbb{E}\left[m(X,{\theta},t^*(\theta))^2\right]-\left[\mathbb{E}\left[m(X,{\theta},
   t^*(\theta))\right]\right]^2.$$
 \end{enumerate}
 \end{theorem}

\begin{remark}{\rm \label{remarque 1}
Let $\theta$ be a given value in $\Theta$. Consider the test
 of the null hypothesis
 \begin{equation}\label{problem de test 1}
 \mathcal{H}_0 : P_0\in\mathcal{M}_\theta \quad
 \text{against}\quad \mathcal{H}_1 : P_0\notin \mathcal{M}_\theta.
 \end{equation}
In view of Theorem \ref{theoreme 3} part 2,  we reject
$\mathcal{H}_0$ against $\mathcal{H}_1$, at asymptotic level
$\alpha\in ]0,1[$, when $2n\widehat{D}_{\varphi
}\left(\mathcal{M}_{\theta },P_{0}\right)$ exceeds the
$(1-\alpha)$- quantile of the $\chi^{2}(l)$ distribution. Theorem
\ref{theoreme 5} part 2 is useful to give an approximation to the
power function
$$P_0\notin \mathcal{M}_\theta \mapsto \beta(P_0):=P_0\left[2n\widehat{D}_{\varphi }\left(\mathcal{M}_{\theta
},P_{0}\right)>q_{(1-\alpha)}\right].$$ We obtain then the
following approximation
 \begin{equation}\label{power approxi 1}
 \beta(P_0)\approx
1-F_\mathcal{N}\left(\frac{\sqrt{n}}{\sigma(\theta)}\left[\frac{q_{1-\alpha}}{2n}-
D_\varphi(\mathcal{M}_\theta,P_0)\right]\right),
\end{equation}
where $F_\mathcal{N}$ is the cumulative distribution function of
the standard normal distribution. From this approximation, we can
give the approximate sample size that ensures a desired power
$\beta$ for a given alternative $P_0\notin \mathcal{M}_\theta$.
Let $n_{0}$ be the positive root of the equation
\[
\beta =1-F_{\mathcal{N}}\left[ \frac{\sqrt{n}}{\sigma
\left(\theta\right) }\left( \frac{q_{(1-\alpha) }}{2n}-D_{\varphi
}\left( \mathcal{M}_\theta,P_0\right) \right) \right]
\]
i.e.,
\[
n_{0}=\frac{\left( a+b\right) -\sqrt{a\left( a+2b\right)
}}{2D_{\varphi }\left(\mathcal{M}_\theta,P_0\right)^{2}}
\]
with $a:=\sigma(\theta^*)^{2}\left[
F_{\mathcal{N}}^{-1}\left(1-\beta\right)\right]^2$ and
$b:=q_{(1-\alpha) } D_{\varphi }\left(
\mathcal{M}_\theta,P_0\right).$ The required sample size is then
$\left\lfloor n_{0}\right\rfloor +1$, where $\left\lfloor
n_{0}\right\rfloor$ denotes the integer part of $n_{0}.$}\\
\end{remark}

\begin{remark}{\rm (\textbf{Generalized empirical likelihood ratio
test}). For testing $\mathcal{H}_0:P_0\in\mathcal{M}_\theta$
against the alternative
$\mathcal{H}_1':\mathcal{M}\setminus\mathcal{M}_\theta$, we
propose to use the statistics
\begin{equation}\label{S n phi}
2nS_{n}^{\varphi}:=2n\left[\widehat{D}_{\varphi }\left(
\mathcal{M}_{\theta },P_{0}\right)
-\inf_{\theta\in\Theta}\widehat{D}_{\varphi }\left(
\mathcal{M}_{\theta },P_{0}\right)\right],
\end{equation}
which converge in distribution  to a $\chi^{2}(d)$ random variable
under $\mathcal{H}_0$ when Assumptions \ref{assumption 1} and
\ref{assumption 2} hold. This can be proved using similar
arguments as in Theorems \ref{theoreme 2} and \ref{theoreme 3}. We
then reject $\mathcal{H}_0$ at asymptotic level $\alpha$ when
$2nS_{n}^{\varphi}>q_{(1-\alpha)}$, the $(1-\alpha)$-quantile of
the $\chi^2(d)$-distribution. Under $\mathcal{H}_1'$ and when
Assumptions \ref{assumption 1},\ref{assumption 2},\ref{assumption
4} and \ref{assumption 5} hold, as in Theorem \ref{theoreme 5}, it
can be proved that
\begin{equation}\label{S n phi alternative}
\sqrt{n}\left( S_{n}^{\varphi }-D_{\varphi }\left(
\mathcal{M}_{\theta },P_{0}\right) \right)
\end{equation}
converges to a centered normal random variable with variance
$$\sigma^{2}(\theta):=\mathbb{E}\left(m(X,\theta ,t^*(\theta))^2\right)-
\left(\mathbb{E}m(X,\theta,t^*(\theta))\right)^2.$$ So, as in the
above Remark, we obtain the following approximation
 \begin{equation}\label{power approxi 2}
 \beta(P_0)\approx
1-F_\mathcal{N}\left(\frac{\sqrt{n}}{\sigma(\theta)}\left[\frac{q_{1-\alpha}}{2n}-
D_\varphi(\mathcal{M}_\theta,P_0)\right]\right)
\end{equation}
to the power function $P_0\in\mathcal{M}/\mathcal{M}_\theta\mapsto
\beta(P_0):=P_0\left[2nS_n^\varphi>q_{(1-\alpha)}\right].$ The
approximated sample size required to achieve a desired power for a
given alternative can be obtained in a similar way.}
\end{remark}

\begin{remark}{\rm (\textbf{Confidence region for the parameter}).
 For a fixed level $\alpha\in]0,1[$,
using convergence (\ref{S n phi}), the set
$$\left\{\theta\in\Theta\text{ such that } 2nS_n^\varphi\leq q_{(1-\alpha)}\right\}$$
is an asymptotic confidence region for $\theta_0$ where
$q_{(1-\alpha)}$ is the $(1-\alpha)$-quantile of the
$\chi^2(d)$-distribution. It is straightforward to see that the
confidence region obtained for the $KL_m$-divergence coincides
with that of \cite{Owen1991} and \cite{Qin-Lawless1994}.}
\end{remark}

\subsection{Asymptotic properties under misspecification}
We address Problem 1 stating the limiting distribution of the
proposed test statistics  under the alternative $\mathcal{H}_1:
P_{0}\notin \mathcal{M}.$ This needs the introduction of
$Q_{\theta^*}^*$, the projection of $P_{0}$ on $\mathcal{M}$.
Assumption 6 below ensures the existence of the ``pseudo-true''
value $\theta^*$ as well as the existence of the projection
$Q_{\theta^*}^*$ of $P_0$ on $\mathcal{M}$, and states some
necessary other regularity conditions. Proposition
\ref{proposition 2} above states the existence and
characterization of the projection $Q_\theta^*$ of $P_0$ on
$\mathcal{M}_\theta$, for a given $\theta\in\Theta$.

\begin{assumption}
\label{assumption 6}  a) $\Theta$ is compact,
$\theta^*:=\arg\inf_{\theta\in\Theta}\sup_{t\in\Lambda_\theta}
\mathbb{E}\left[m(X,\theta,t)\right]$ exists and is unique; b)
$g(X,\theta)$ is continuous  at each $\theta\in\Theta$ with
probability one;\\ c) $\mathbb{E}\left[\sup_{\{\theta\in\Theta,
t\in N_{t^*(\theta)}\}}|m(X,\theta,t)|\right]<\infty,$ where
$N_{t^*(\theta)}\subset \Lambda_\theta$  is a compact set such
that $t^*(\theta)\in\text{int}\left(N_{t^*(\theta)}\right)$; d)
for all $\theta\in\Theta$, the functions
$\mathds{1}_{\mathbb{R}^m}, g_1,\ldots,g_l$  are linearly
independent in the  sense that
$P_0\left\{x\in\mathbb{R}^m~|~t_0+\sum_{j=1}^lt_jg_j(x,\theta)
\neq 0\right\}>0$,  for all $t\in\mathbb{R}^{1+l}$ with $t\neq 0.$
\end{assumption}

\begin{remark}\label{remark 2 sur les conditions}
 Assumption \ref{assumption 6}.d ensures the strict concavity of the function
 $t\in\Lambda_\theta \mapsto \mathbb{E}\left[m(X,\theta,t)\right]$
 on the convex set $\Lambda_\theta$, which implies the uniqueness
 of $t^*(\theta)$, for all $\theta\in\Theta$. This Assumption can
 be replaced by the following one : for all $\theta\in\Theta$, there
 exists a neighborhood $N_{t^*(\theta)}$ of $t^*(\theta)$ such that
 $$\mathbb{E}[\sup_{t\in N_{t^*(\theta)}} \left\|\partial m(X,\theta,t)/\partial
 t\right\|]<\infty, \quad\mathbb{E} [\sup_{t\in N_{t^*(\theta)}} \left\|\partial^2
m(X,\theta,t)/\partial
 t^2\right\|]<\infty$$ and the matrix  $\mathbb{E}\left[\partial^2 m(X,\theta,t^*(\theta))/\partial
 t^2\right]<\infty$ is nonsingular, which implies the uniqueness
 of $t^*(\theta)$, for all $\theta\in\Theta$.
\end{remark}

\begin{theorem}\label{theoreme 6}
 Under Assumption \ref{assumption 6}, we have
 \begin{enumerate}
  \item [1)] $\|\widehat{t}(\theta)-t^*(\theta)\|$ converges in
  probability to $0$ uniformly in $\theta\in\Theta$.
  \item [2)] $\widehat{\theta}_\varphi$ converges in probability to
  $\theta^*$;
  \item [3)] $\widehat{D}_\varphi(\mathcal{M},P_0)$ converges in
  probability to $D_\varphi(\mathcal{M},P_0)$.
 \end{enumerate}
\end{theorem}

\noindent The asymptotic normality of the test  statistics under
misspecification requires the following additional  conditions.

\begin{assumption}
  \label{assumption 7}
 a) $\theta ^{\ast }\in \text{int}(\Theta )$;
b) there exists $\mathcal{N}\subset \Theta\times\Lambda_\Theta$,
some compact  neighborhood of
 $(\theta^*,t^*(\theta^*))$, such that with probability one  $(\theta, t)\in \mathcal{N}\mapsto m(X,\theta,t)$
 is $\mathcal{C}^2$ and
  $$\mathbb{E}[ \sup_{(\theta,t)\in\mathcal{N}}\|\partial m(X,\theta,t)/\partial (\theta, t)\|]<\infty, \quad
  \mathbb{E}[ \sup_{(\theta,t)\in\mathcal{N}}\|\partial^2m(X,\theta,t)/\partial (\theta, t)^2\|]<\infty;$$
    c) as $\delta\to 0$,
    $$\mathbb{E}\left\{\sup_{\{(t,\theta); \|(t,\theta)-(t^*(\theta^*),\theta^*)\|\leq
    \delta\}}\left\|\partial^2m(X,\theta,t)/\partial (\theta,t)^2-
    \partial^2m(X,\theta^*,t^*(\theta^*))/\partial (\theta,t)^2\right\|\right\}\to 0;$$
 d) $\mathbb{E}\left[ m(X,\theta ^{\ast },t^{\ast
}(\theta ^{\ast }))^{2}\right] ,$ $\mathbb{E}\left[ \left\Vert
\partial m(X,\theta ^{\ast },t^{\ast }(\theta ^{\ast }))/\partial
t\right\Vert ^{2}\right] $ and $\mathbb{E}\left[ \left\Vert
\partial m(X,\theta ^{\ast },t^{\ast }(\theta ^{\ast })/\partial
\theta \right\Vert ^{2}\right] $ are finite, and the matrix
\[
S:=\left(
\begin{array}{cc}
S_{11} & S_{12} \\
S_{21} & S_{22} \\
\end{array}
\right) ,
\]
is  nonsingular, where $S_{11}:=\mathbb{E}\left[
\partial ^{2}m(X,\theta ^{\ast },t^{\ast }(\theta ^{\ast
}))/\partial t^{2}\right] $,\\
$S_{12}={S_{21}}^{\top}:=\mathbb{E}\left[ \partial ^{2}m(X,\theta
^{\ast },t^{\ast }(\theta ^{\ast }))/\partial t\partial \theta
\right] $ and $S_{22}:=\mathbb{E}\left[ \partial ^{2}m(X,\theta
^{\ast },t^{\ast }(\theta ^{\ast }))/\partial \theta ^{2}\right].$
\end{assumption}

\begin{remark}\label{remark 2 sur les conditions de 3 ordre}
 Assumption \ref{assumption 7}.c is used here to relax the
 condition on the third derivatives (in $t$ and $\theta$) of the function $(\theta,t)\mapsto
 m(X,\theta,t)$.
\end{remark}

\begin{theorem}\label{theoreme 7}
  Under Assumptions \ref{assumption 6} and \ref{assumption 7}, we
  have
 \begin{enumerate}
  \item [1)]
  \begin{equation*}
    \sqrt{n}\left(
\begin{array}{c}
  \widehat{t}(\widehat{\theta}_\varphi)-t^*(\theta^*) \\
  \widehat{\theta}_\varphi-\theta^* \\
\end{array}
\right)
\end{equation*}
converges in distribution to a centered  normal random vector with
covariance matrix
\begin{equation*}\label{matrice limit a l'exter du model}
  W := S^{-1}M{S^{-1}}
\end{equation*}
where
\begin{equation*}
M := \mathbb{E}\left[\left[
\begin{array}{c}
  \frac{\partial}{\partial t}m\left(X,\theta^*,t^*(\theta^*)\right) \\
  \frac{\partial }{\partial \theta}m\left(X,\theta^*,t^*(\theta^*)\right) \\
\end{array}
\right]\left[
\begin{array}{c}
  \frac{\partial}{\partial t}m\left(X,\theta^*,t^*(\theta^*)\right) \\
  \frac{\partial }{\partial \theta}m\left(X,\theta^*,t^*(\theta^*)\right) \\
\end{array}
\right]^{\top}\right];
\end{equation*}
\item [2)]
$\sqrt{n}\left(\widehat{D}_\varphi(\mathcal{M},P_0)-D_\varphi(\mathcal{M},P_0)\right)$
converges in distribution to a centered normal random variable
with variance
$$\sigma^2(\theta^*):=\mathbb{E}\left[m(X,\theta^*,t^*(\theta^*))^2\right]-
\left[\mathbb{E}\left[m(X,\theta^*,t^*(\theta^*))\right]\right]^2.$$
\end{enumerate}
\end{theorem}

\begin{remark}{\rm
 In the case of EL, i.e., when $\varphi(x)=-\log x +x -1$,
 Assumption \ref{assumption 6}.c implies that
 \begin{equation}\label{bounded}
 -\infty < \inf_{x\in\mathbb{R}^m} t_0+t^{{\top}} g(x,\theta)\leq
 \sup_{x\in\mathbb{R}^m}t_0+t^{\top}g(x,\theta)<1
 \end{equation}
 for all $x\in\mathbb{R}^m-P_0$-a.s.,  for all $\theta\in \Theta$ and for all $t\in N_{t^*(\theta)}$. This
 imposes a restriction on the model when the support of $P_0$ and the functions $g_j$
 are unbounded. Indeed, when the support of $P_0$
 is for example the  whole space $\mathbb{R}^m$, the condition above
 does not hold when $g$ is unbounded. In this case, the EL estimate may cease to be consistent
 as it is stated by \cite{Schennach2007} under misspecification. This is a potential problem
 for all divergences associated to $\varphi$-functions with domain
 of the form $(a,+\infty[$, $]-\infty,b)$ or $(a,b),$ where $a$ and
 $b$ are some finite real numbers; it is the case of  modified $\chi^2$,
 Hellinger,  KL and modified $KL$ divergences.
 At the contrary, Assumption 6.c may be satisfied for other divergences
 associated to $\varphi$ functions with
$\text{dom}
 \varphi=\mathbb{R}$ which
 is the case of $\chi^2$ divergence for example.}
\end{remark}

\begin{remark} \label{remark power approximation}
{\rm Theorem \ref{theoreme 7} part 2 is useful for the computation
of the power function. For  testing the null hypothesis
$\mathcal{H}_0 : P_0\in\mathcal{M}$ against the alternative
$\mathcal{H}_1:P_0\notin\mathcal{M}$, the power function is
\begin{equation}\label{power function model}
P_0\notin\mathcal{M}\mapsto
\beta(P_0):=P_0\left[2n\widehat{D}_{\varphi }\left(
\mathcal{M},P_0\right)
>q_{(1-\alpha)}\right].
\end{equation}
Using Theorem \ref{theoreme 7} part 2, we obtain the following
approximation to the power function (\ref{power function model}):
\begin{equation}\label{power approxi}
\beta(P_0)\thickapprox 1-F_{\mathcal{N}}\left[
\frac{\sqrt{n}}{\sigma\left( \theta ^{\ast }\right) }\left(
\frac{q_{(1-\alpha)}}{2n}- D_\varphi\left( \mathcal{M},P_0\right)
\right) \right]
\end{equation}
where $F_\mathcal{N}$ is the empirical cumulative distribution of
the standard normal distribution. From the proxy value of
$\beta(P_0)$ hereabove,   the approximate sample size that ensures
a given power $\beta$ for a given alternative $P_0\not\in
\mathcal{M}$ can be obtained as follows. Let $n_{0}$ be the
positive root of the equation
\[
\beta =1-F_{\mathcal{N}}\left[ \frac{\sqrt{n}}{\sigma( \theta^*)
}\left( \frac{q_{(1-\alpha) }}{2n}-D_{\varphi }\left(
\mathcal{M},P_0\right) \right) \right]
\]
i.e.,
\[
n_{0}=\frac{\left( a+b\right) -\sqrt{a\left( a+2b\right)
}}{2D_\varphi\left( \mathcal{M},P_0\right) ^{2}},
\]
where  $a:=\sigma(\theta^*) ^{2}\left[ F_{\mathcal{N}}^{-1}\left(
1-\beta \right) \right] ^{2}$ and $b:=q_{(1-\alpha) } D_{\varphi
}\left( \mathcal{M},P_0\right).$ The required sample size is then
$\left\lfloor n_{0}\right\rfloor +1$.}
\end{remark}

\section{Simulation results: Approximation of the power function of
the empirical likelihood ratio test}
\begin{figure}[!h]
 \caption[]{Approximation of the power function}
\centerline{
  \begin{tabular}{ c  c }
 \includegraphics[width=.55\textwidth]{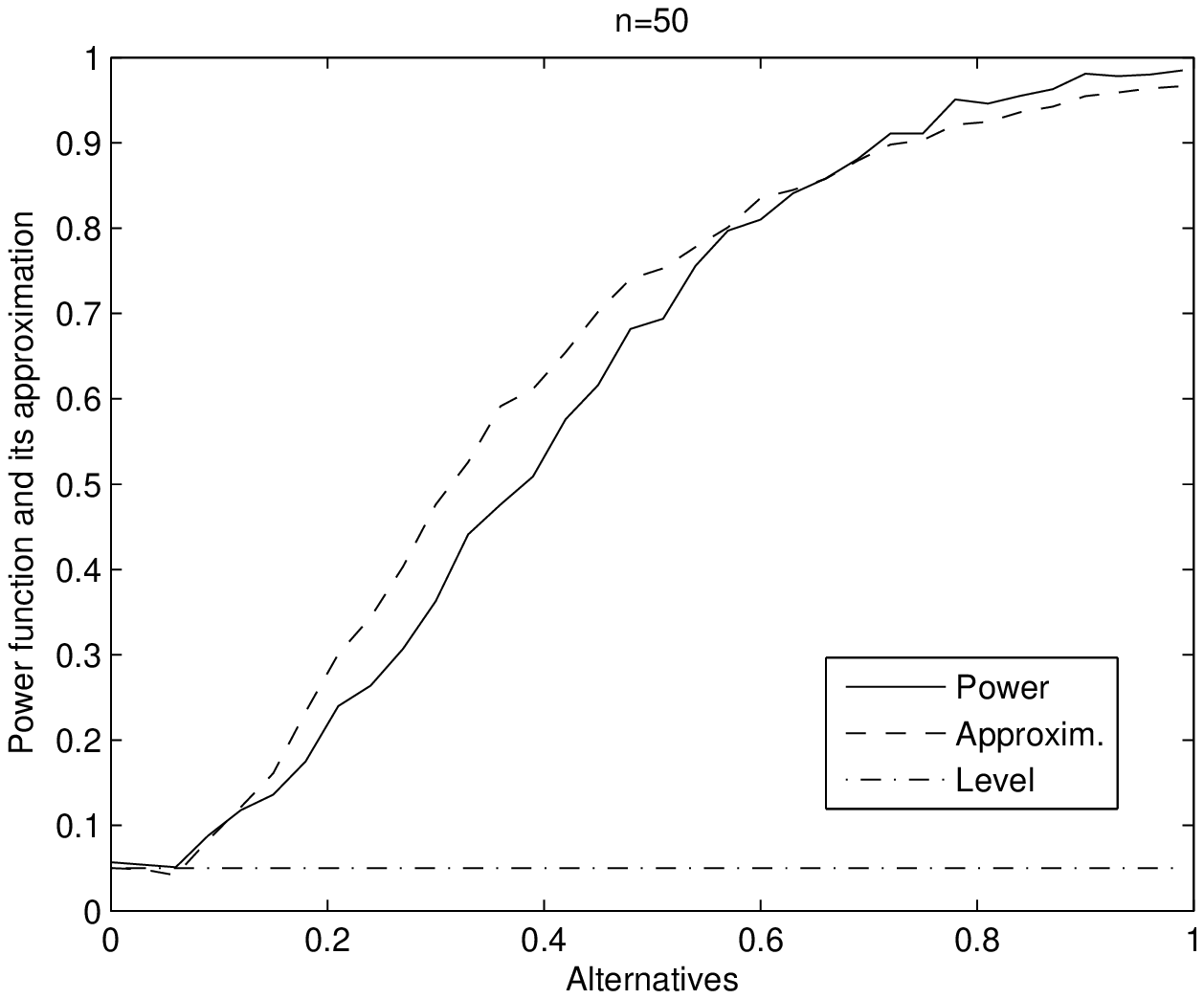}
    &   \includegraphics[width=.55\textwidth]{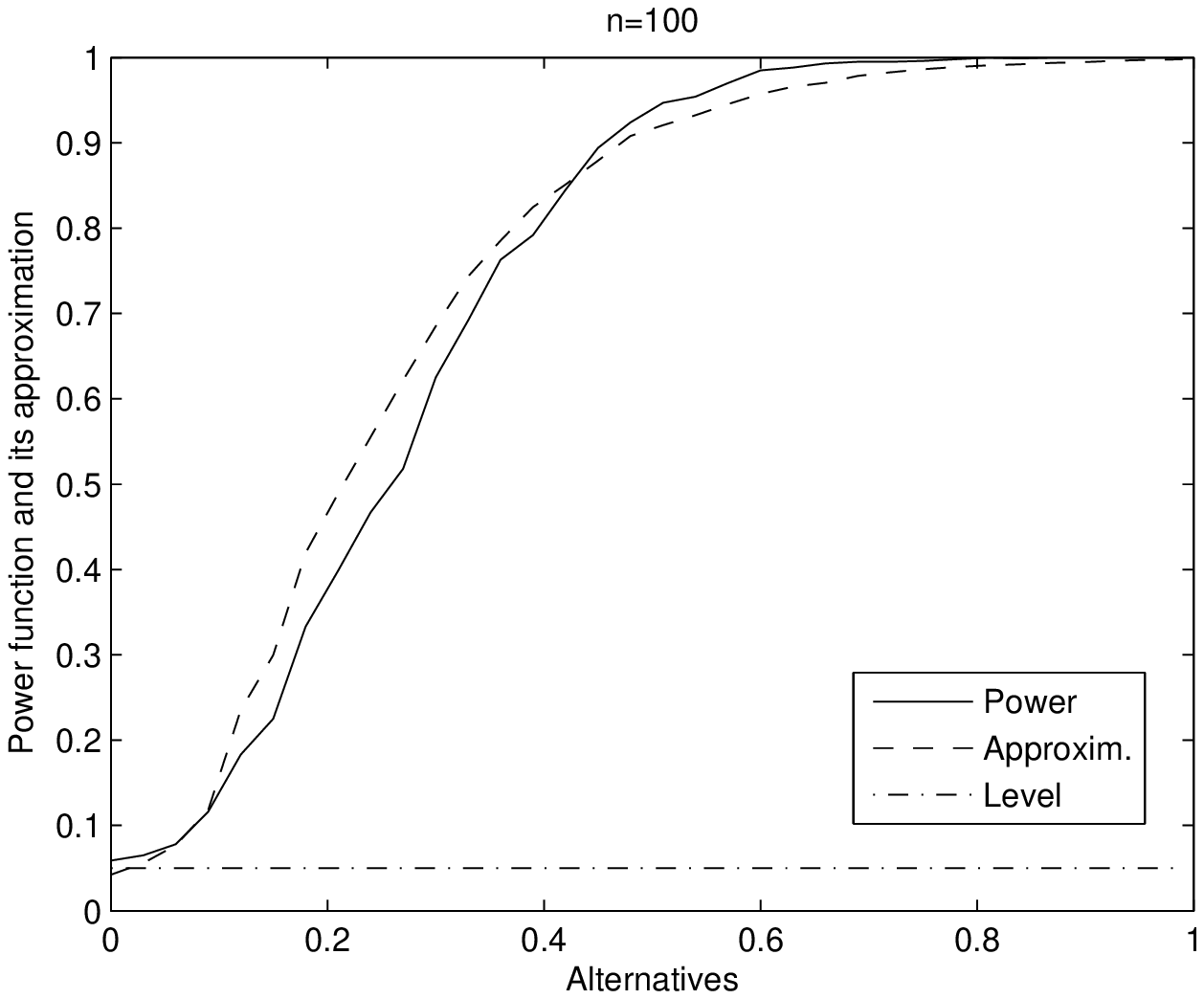}
   \end{tabular} }
\par
\centerline{
  \begin{tabular}{ c  c }
 \includegraphics[width=.55\textwidth]{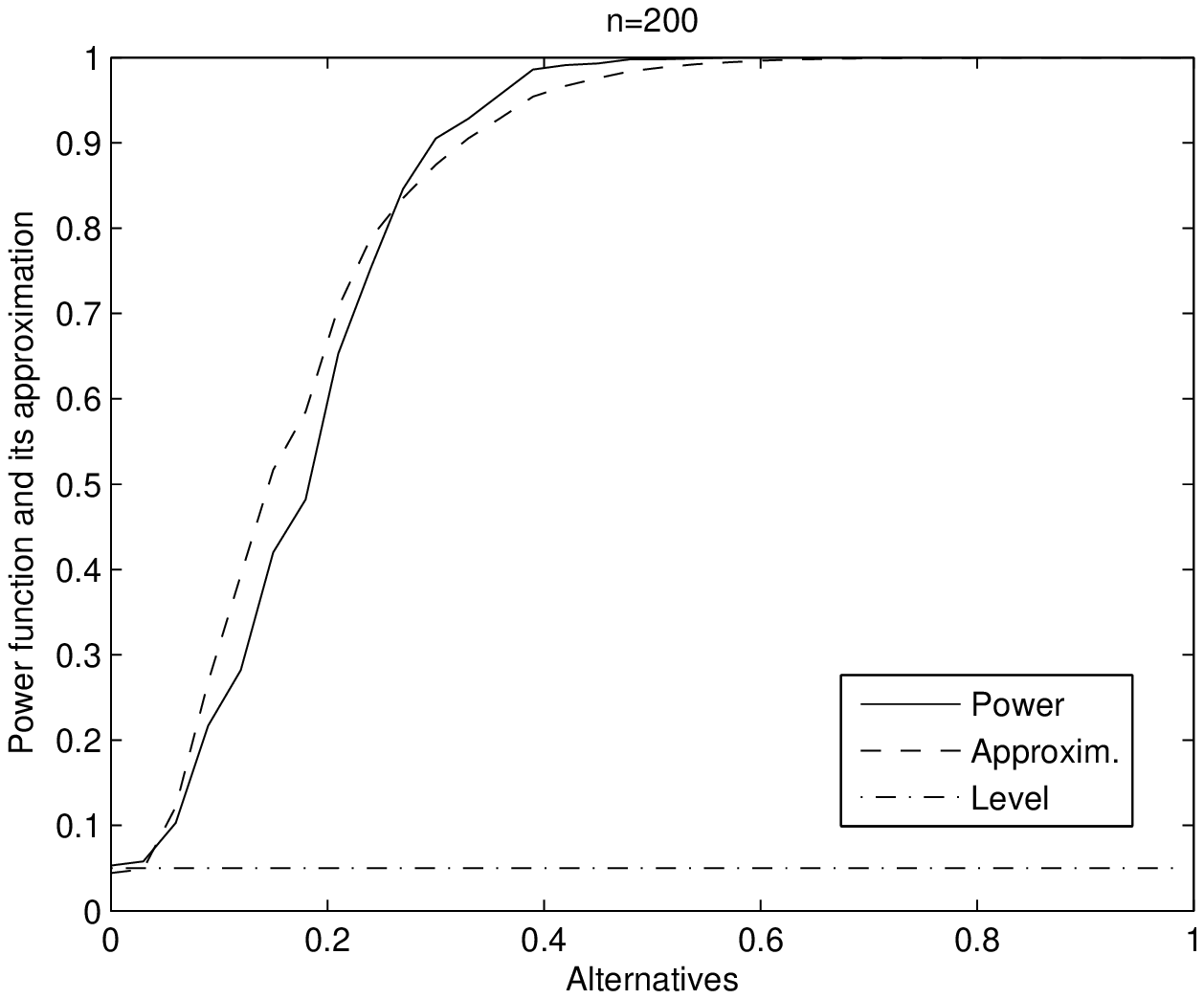}
    &   \includegraphics[width=.55\textwidth]{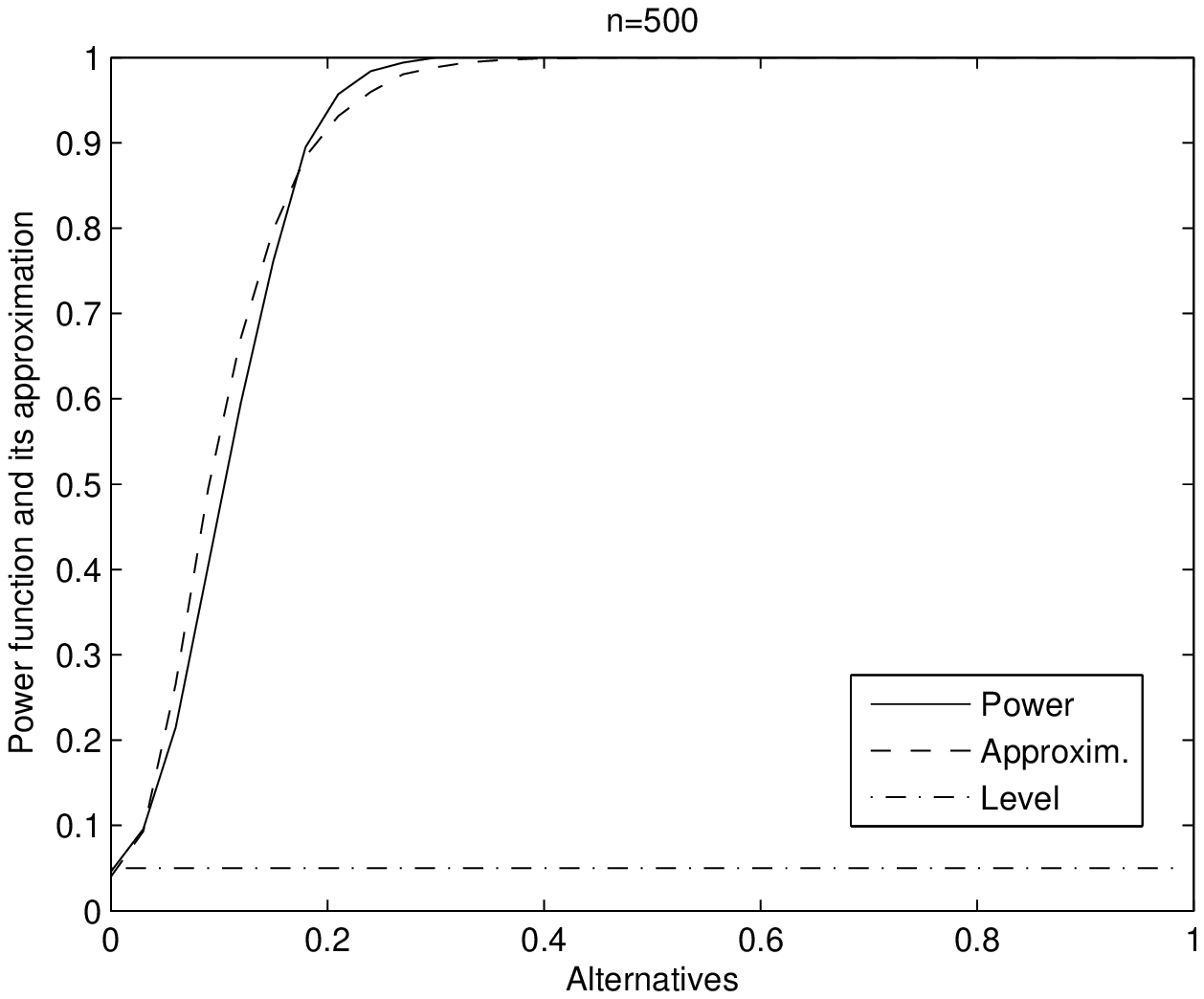}
   \end{tabular} }
 \label{fig power approxim}
\end{figure}

\noindent  We will illustrate by simulation the accuracy of the
power
 approximation (\ref{power approxi}) in the case of EL method, i.e., when
 $\varphi(x)=-\log x +x-1.$
 Consider the test problem of the composite null hypothesis
 $$\mathcal{H}_0: P_0\in\mathcal{M}\quad \text{against the alternative}
 \quad \mathcal{H}_1:P_0\notin \mathcal{M},$$
 where
$\mathcal{M}:=\bigcup_{\theta\in\mathbb{R}}\mathcal{M}_\theta$ and
 $\mathcal{M}_\theta$ is the set of all s.f.m's satisfying the constraints  $\int dQ(x)=1$
 and $\int g(x,\theta)~dQ(x)=0$ with $g(x,\theta):=(x,
 x^2-\theta)^{\top}$, namely
  $$\mathcal{M}_\theta:=\left\{Q\in M ~ \text{ such that } ~ \int_\mathbb{R} dQ(x)=1 \text{ and }
  \int_\mathbb{R} g(x,\theta)~dQ(x)=0  \right\},$$
  where $\theta\in\mathbb{R}$ is the parameter of interest. We
  consider the asymptotic level $\alpha=0.05$ and the alternatives
  $P_0:=\mathcal{U}([-1,1+\epsilon])\not\in\mathcal{M}$ for different values of
  $\epsilon$ in the interval $]0,1]$. Note that when $\epsilon=0$ then the
  uniform distribution $\mathcal{U}([-1,1])$ belongs to
  the model $\mathcal{M}$. For this model, we can show also that all Assumptions  of Theorem \ref{theoreme 2}
are satisfied when $\epsilon=0$, and all Assumptions of Theorem
\ref{theoreme 7}  are met under alternatives. In Figure
 \ref{fig power approxim}, the power function (\ref{power function
model}) is plotted (with a continuous line),
  with sample sizes $n=50, n=100$, $n=200$ and $n=500$, for different values
  of $\epsilon$. Each power entry was obtained by Monte-Carlo from
  $1000$ independent runs. The approximation (\ref{power approxi}) is
  plotted (with a dashed line) as a function of $\epsilon$.
 The estimates $\widehat{\theta}_\varphi$  and  $\widehat{D}_\varphi(\mathcal{M},P_0)$
 are calculated using the Newton-Raphson algorithm. We observe from
Figure \ref{fig power approxim} that the approximation is accurate
even for moderate sample sizes.

\section{Concluding remarks and possible developments}
\noindent We have proposed new estimates and tests for model
satisfying linear constraints with unknown parameter through
divergence based methods which generalize the EL approach. This
leads to the obtaining of the limiting distributions of the test
statistics and the estimates  under alternatives and under
misspecification. Consistency of the test statistics under the
alternatives is the starting point for the study of the optimality
of the tests through Bahadur approach; also the generalized
Neyman-Pearson optimality of EL test (as developed by
\cite{Kitamura2001}) can be adapted for empirical divergence based
methods. Many problems remain to be studied in the future such as
the choice of the divergence which leads to an optimal (in some
sense) estimator or test in terms of efficiency and/or robustness.
Preliminary simulation results show that Hellinger divergence
enjoys good properties in terms of efficiency-robustness; see
\cite{BK_2006_div_emp}. Also comparisons under local alternatives
should be developed.

\section{Appendix}
\noindent \textbf{Proof of Theorem \ref{theoreme 1}}.\\ The same
arguments, used for the proof of Theorem 3.1 in
\cite{NeweySmith2004}, hold when their criterion function
$(\theta,\lambda)\in\Theta\times\mathbb{R}^l\mapsto
\frac{1}{n}\sum_{i=1}^n \rho(\lambda^{\top}g(X,\theta))$ is
replaced by our function
$(\theta,t)\in\Theta\times\mathbb{R}^{1+l}\mapsto
\frac{1}{n}\sum_{i=1}^n m(t^{\top}\overline{g}(X,\theta)).$ In
particular, we have $$\max_{1\leq i\leq
n}\left|\widehat{t}(\widehat{\theta}_\varphi)^{\top}\overline{g}(X_i,\widehat{\theta}_\varphi)\right|\to
0$$ in probability, which implies that
$\widehat{t}(\widehat{\theta}_\varphi)\in\text{int}(\Lambda^{(n)}_{\widehat{\theta}_\varphi})$
with probability one as $n\to\infty$, since $a^*<0<b^*$.\\

\noindent \textbf{Proof of Theorem \ref{theoreme 2}}.\\ The proof
is similar to that of \cite{NeweySmith2004} Theorem 3.2. Hence, it
is
omitted.\\

\noindent \textbf{Proof of Theorem \ref{theoreme 3}}.\\
  It is a particular case of Theorem \ref{theoreme 1} taking
$\Theta=\{\theta\}.$ Hence, the proof is omitted.\\

\noindent \textbf{Proof of Theorem \ref{theoreme 4}}.\\ 1) First,
note that $t^*(\theta)$ exists and is unique by Assumption
\ref{assumption 4}. By the uniform weal law of large numbers
(UWLLN), using continuity of $m(X,\theta,t)$ in $t$, and
Assumption \ref{assumption 4}.b, we obtain
\begin{equation}\label{eqn3}
\left|P_nm(\theta,t)-\mathbb{E}\left[m(X,\theta,t)\right]\right|\to
0,\end{equation} in probability uniformly in $t$ over the compact
set $N_{t^*(\theta)}$. Using this and the fact that
$t^*(\theta):=\arg\sup_{t\in\Lambda_\theta}P_0m(\theta,t)$ is
unique and belongs to $\text{int}(N_{t^*(\theta)})$ and the strict
concavity of $t\mapsto P_0m(\theta,t)$, we conclude that  any
value
 \begin{equation}\label{eqn t bare}\overline{t}:=\arg\sup_{t\in N_{t^*(\theta)}}P_nm(\theta,t)
 \end{equation} converges in
probability to $t^*(\theta)$; see e.g. Theorem 5.7 in
\cite{vanderVaart1998}. We end then the proof by showing that
$\widehat{t}(\theta)$ belongs to $\text{int}(N_{t^*(\theta)})$
with probability one as $n\to \infty$, and therefore it converges
to $t^*(\theta).$ In fact, since for $n$ sufficiently large any
value $\overline{t}$ lies in the interior of $N_{t^*(\theta)}$,
concavity of $t\mapsto P_nm(\theta,t)$ implies that no other point
$t$ in the complement of $\text{int}(N_{t^*(\theta)})$ can
maximize $P_nm(\theta,t)$ over $t\in\mathbb{R}^{1+l}$, hence
$\widehat{t}(\theta)$ must belongs
to $\text{int}(N_{t^*(\theta)})$.\\
 2) With probability tending to $1$ as $n\to \infty$, we have
$\widehat{D}_\varphi(\mathcal{M}_{\theta},P_0)=P_nm(\theta,\widehat{t})=P_nm(\theta,\overline{t})$.
Hence, we can write
$$\left|\widehat{D}_\varphi(\mathcal{M}_{\theta},P_0)-D_\varphi(\mathcal{M}_{\theta},P_0)\right|=
\left|P_nm(\theta,\overline{t})-P_0m(\theta,t^*(\theta))\right|=:|A|,$$
and
$$P_nm(\theta,t^*(\theta))-P_0m(\theta,t^*(\theta))\leq A\leq
P_nm(\theta,\overline{t})-P_0m(\theta,\overline{t}).$$
Both the RHS and the LHS in the above display tend to $0$ in
probability by (\ref{eqn3}). Hence,
$\left|\widehat{D}_\varphi(\mathcal{M}_{\theta},P_0)-D_\varphi(\mathcal{M}_{\theta},P_0)\right|$
 tends to $0$ in probability as $n\to\infty$. This ends the proof.\\

\noindent \textbf{Proof of Theorem \ref{theoreme 5}}.\\ 1) For $n$
sufficiently large, by a Taylor expansion, there exists
$\overline{t}\in \mathbb{R}^{1+l}$ inside the segment that links
$\widehat{t}$ and $t^*(\theta)$ with
\begin{equation}\label{Taylor Expansion}
    \begin{array}{ccl}
    0 & = & P_nm'(\theta,\widehat{t})\\
      & = &
    P_nm'(\theta,t^*(\theta))+\left(P_n m''(\theta,\overline{t})\right)^{\top}
    \left(\widehat{t}-t^*(\theta)\right).
    \end{array}
\end{equation}
By Assumptions \ref{assumption 5}.a and \ref{assumption 5}.b,
using the fact that $\overline{t}=t^*(\theta)+o_P(1)$ and the
UWLLN, we can prove that $$P_n m''(\theta,\overline{t})= P_0
m''(\theta,t^*(\theta))+o_P(1).$$ Using this display, one gets
from (\ref{Taylor Expansion})
\begin{equation}\label{equation loi limite}
    -P_nm'(\theta,t^*(\theta))=\left(P_0m''(\theta,t^*(\theta))+o_P(1)\right)
    \left(\widehat{t}-t^*(\theta)\right).
\end{equation}
Assumptions \ref{assumption 4}.a and \ref{assumption 5}.a imply
that $P_0m'(\theta,t^*(\theta))=0$. Hence, by the central limit
theorem (CLT), we have
$$\sqrt{n}P_nm'(\theta,t^*(\theta))=O_P(1),$$
which by (\ref{equation loi limite}) implies that
$\sqrt{n}\left(\widehat{t}-t^*(\theta)\right)=O_P(1).$ Hence, from
(\ref{equation loi limite}), we get
\begin{equation}\label{equation qui donne la loi limite}
    \sqrt{n}\left(\widehat{t}-t^*(\theta)\right)=
    {\left[-P_0m''(\theta,t^*(\theta))\right]}^{-1}
    \sqrt{n}P_nm'(\theta,t^*(\theta))+o_P(1).
\end{equation}
The CL and Slutsky theorems conclude the proof of part 1.\\
 2)  Using the fact that
$\left(\widehat{t}-t^*(\theta)\right)=O_P(1/\sqrt{n})$ and
$P_nm'(\theta,t^*(\theta))=P_0m'(\theta,t^*(\theta))+o_P(1)=0+o_P(1)=o_P(1)$,
 we obtain
\begin{eqnarray}
\sqrt{n}\left(\widehat{D}_\varphi(\mathcal{M}_\theta,P_0)-
D_\varphi(\mathcal{M}_\theta,P_0)\right) & = & \sqrt{n}\left(
P_nm(\theta,\widehat{t})-P_0m(\theta,t^*(\theta))\right)\nonumber\\
 & = &
 \sqrt{n}\left(P_{n}m(\theta,t^*(\theta))-P_0m(\theta,t^*(\theta))\right)+o_{P}(1),\nonumber
\end{eqnarray}
and the CL and Slutsky theorems conclude the proof.\\

\noindent \textbf{Proof of Theorem \ref{theoreme 6}}.\\ 1) First
note that Assumption \ref{assumption 6}.d implies that the
function $t\in\Lambda_\theta\mapsto \mathbb{E} m(X,\theta,t)$ is
strictly concave for all $\theta\in\Theta$, which implies that
$t^*(\theta)$ is unique for all $\theta\in\Theta$. By the UWLLN,
using continuity of $m(X,\theta,t)$, in $\theta$ and $t$, and
Assumption \ref{assumption 6}.c, we obtain the uniform convergence
in probability, over the compact set
$\left\{(\theta,t)\in\Theta\times\mathbb{R}^{1+l};
~\theta\in\Theta, t\in N_{t^*(\theta)}\right\}$,
 \begin{equation}\label{eqn 4}\sup_{\{\theta\in\Theta,t\in N_{t^*(\theta)}\}}\left|P_nm(\theta,t)
-P_0m(\theta,t)\right|\to 0.
  \end{equation} We can then prove
the convergence in probability
$\sup_{\theta\in\Theta}\|\widehat{t}(\theta)-t^*(\theta)\|\to 0$
in two steps. Step 1: Let $\eta >0$. We will show that
$P_0\left[\sup_{\theta\in\Theta}\|\overline{t}(\theta)-t^*(\theta)\|\geq
\eta\right]\to 0$ for any value
 \begin{equation}\label{eqn t bare theta}
 \overline{t}(\theta):=\arg\sup_{t\in N_{t^*(\theta)}}P_nm(\theta,t).
 \end{equation}
Step 2: To conclude the proof, we will show that
$\widehat{t}(\theta)$ belongs to $\text{int}(N_{t^*(\theta)})$
with probability one as $n\to\infty$  for all $\theta\in\Theta$.
Let $\eta
>0$ such that
$\sup_{\theta\in\Theta}\|\overline{t}(\theta)-t^*(\theta)\|\geq
\eta$. Sine $\Theta$ is a compact set, by continuity there exists
$\overline{\theta}\in\Theta$ such that
$\sup_{\theta\in\Theta}\|\overline{t}(\theta)-t^*(\theta)\|=
\|\overline{t}(\overline{\theta})-t^*(\overline{\theta})\|\geq
\eta$. Hence, there exists $\varepsilon >0$  such that
$P_0m(\overline{\theta},t^*(\overline{\theta}))-
P_0m(\overline{\theta},\overline{t}(\overline{\theta}))>\varepsilon$.
In fact, $\varepsilon$ may be defined as follows
$$\varepsilon := \inf_{\theta\in\Theta}
\sup_{\{t\in N_{t^*(\theta)};~ \|t-t^*(\theta)\|\geq \eta\}}
\mathbb{E}[m(X,\theta,t^*(\theta))]- \mathbb{E}[m(X,\theta,t)],$$
which is strictly positive by the strict concavity of
$\mathbb{E}[m(X,\theta,t)]$ in $t$ for all $\theta\in\Theta$, the
uniqueness of $t^*(\theta)\in\text{int}(N_{t^*(\theta)})$ and the
fact that $\Theta$ is compact. Hence the event
$\left[\sup_{\theta\in\Theta}\|\overline{t}(\theta)-t^*(\theta)\|\geq
\eta\right]$ implies the event
$$\left[P_0m(\overline{\theta},t^*(\overline{\theta}))-
P_0m(\overline{\theta},\overline{t}(\theta))\geq \varepsilon
\right],$$ from which we obtain
\begin{equation}\label{eqn 5}
 P_0\left[\sup_{\theta\in\Theta}\|\overline{t}(\theta)-t^*(\theta)\|\geq
\eta\right]\leq P_0\left[
P_0m(\overline{\theta},t^*(\overline{\theta}))-
P_0m(\overline{\theta},\overline{t}(\theta))\geq
\varepsilon\right].
\end{equation}
On the other hand, by (\ref{eqn 4}), we have
\begin{eqnarray} P_0m(\overline{\theta},t^*(\overline{\theta}))-
P_0m(\overline{\theta},\overline{t}(\theta)) & = &
P_nm(\overline{\theta},t^*(\overline{\theta}))-
P_0m(\overline{\theta},\overline{t}(\theta))+o_P(1)\nonumber\\
& \leq & P_nm(\overline{\theta},\overline{t}(\overline{\theta}))-
P_0m(\overline{\theta},\overline{t}(\theta))+o_P(1)\nonumber\\
& \leq & \sup_{\{\theta\in\Theta, t\in N_{t^*(\theta)}\}}
|P_nm(\theta,t)-P_0m(\theta,t)|+o_P(1).\nonumber
\end{eqnarray}
Combining this with (\ref{eqn 5}) and (\ref{eqn 4}), we conclude
that
 \begin{equation}\label{convergence en proba de t bare theta}
\sup_{\theta\in\Theta}\|\overline{t}(\theta)-t^*(\theta)\|\to 0
 \end{equation}
in probability. In particular,
$\overline{t}(\theta)\in\text{int}(N_{t^*(\theta)})$ for
sufficiently large $n$, uniformly in  $\theta\in \Theta$. Since
$t\mapsto P_nm(\theta,t)$ is concave, then the maximizer
$\widehat{t}(\theta)$ belongs to $\text{int}(N_{t^*(\theta)})$ for
sufficiently large $n$; hence the
same result (\ref{convergence en proba de t bare theta}) holds when
$\overline{t}(\theta)$ is replaced by $\widehat{t}(\theta)$.\\
  2) From part 1, we have for  large $n$,
 \begin{eqnarray}
 \sup_{\theta\in\Theta}|P_nm(\theta,\widehat{t}(\theta))-P_0m(\theta,t^*(\theta))|
 & = &
 \sup_{\theta\in\Theta}|P_nm(\theta,\overline{t}(\theta))-P_0m(\theta,t^*(\theta))|=:
 \sup_{\theta\in\Theta}|B|.\nonumber
 \end{eqnarray}
 On the other hand, we have
  $$P_nm(\theta,t^*(\theta))-P_0m(\theta,t^*(\theta))\leq B\leq P_nm(\theta,\overline{t}(\theta))
  -P_0m(\theta,\overline{t}(\theta)).$$
 By Assumption \ref{assumption 6}.c, and the convergence in
 probability $\sup_{\theta\in\Theta}\|\overline{t}(\theta)-t^*(\theta)\|\to
 0$,
 both the  RHS  and LHS of the above display
 tends to $0$ in probability uniformly in $\theta\in\Theta$, by the UWLLN. Hence,
 $\sup_{\theta\in\Theta}|P_nm(\theta,\widehat{t}(\theta))-P_0m(\theta,t^*(\theta))|\to 0$
 in probability. Now, since the minimizer $\theta^*$ of $\theta\mapsto
P_0m(\theta,t^*(\theta))$
 over the compact set $\Theta$
 is unique and interior point of $\Theta$, by continuity and the above uniform convergence, we
 conclude that $\widehat{\theta}_\varphi$ tends  in probability to
 $\theta^*$; see e.g. Theorem 5.7 in
\cite{vanderVaart1998}.\\
   3) This holds as a consequence of the uniform convergence in
  probability
  \begin{equation}\label{conv unif en theta}
  \sup_{\theta\in\Theta}|P_nm(\theta,\widehat{t}(\theta))-P_0m(\theta,t^*(\theta))|\to 0
  \end{equation}
   proved in part 2 above. In fact, we have for $n$
  sufficiently large
    $$ |\widehat{D}_\varphi(\mathcal{M}, P_0)-D_\varphi(\mathcal{M}, P_0)|=
    |P_nm(\widehat{\theta},\widehat{t}(\widehat{\theta}))-P_0m(\theta^*,t^*(\theta^*))|=:|C|,$$
 with
  $$ P_nm(\widehat{\theta},\widehat{t}(\widehat{\theta}))-P_0m(\widehat{\theta},t^*(\widehat{\theta}))
   \leq C
   \leq P_nm(\theta^*,\widehat{t}(\theta^*))-P_0m(\theta^*,t^*(\theta^*))$$
  and both the RHS and LHS tend to $0$ in probability by (\ref{conv unif en
  theta}). This concludes the proof.\\

\noindent \textbf{Proof of Theorem \ref{theoreme 7}}.\\
  1) By the
first order conditions, with probability tending to one, we have
\begin{equation*}
\left\{
\begin{array}{lll}
P_n\frac{\partial}{\partial
t}m\left(\widehat{\theta},\widehat{t}(\widehat{\theta})\right)
& = & 0\\
P_n\frac{\partial}{\partial\theta}
m\left(\widehat{\theta},\widehat{t}(\widehat{\theta})\right)+P_n
\frac{\partial}{\partial t}
m\left(\widehat{\theta},\widehat{t}(\widehat{\theta})\right)
\frac{\partial}{\partial \theta}\widehat{t}(\widehat{\theta}) & =
& 0.
\end{array}
\right.
\end{equation*}
The second term in the LHS of the second equation is equal to $0$,
due to the first equation. Hence, $\widehat{t}(\widehat{\theta})$
and $\widehat{\theta}$ are solutions of the somehow simpler system
\begin{eqnarray}
P_n\frac{\partial}{\partial
t}m\left(\widehat{\theta},\widehat{t}(\widehat{\theta})\right)
& = & 0 \label{E1} \\
P_n\frac{\partial}{\partial\theta}
m\left(\widehat{\theta},\widehat{t}(\widehat{\theta})\right) & = &
 0. \label{E2}
\end{eqnarray}
 Using a Taylor expansion in (\ref{E1}) in
$(\widehat{\theta},\widehat{t})$ around $(\theta^*,t^*)$; there
exists $\left(\overline{\theta},\overline{t}\right)$ inside the
segment that links
$(\widehat{\theta},\widehat{t}(\widehat{\theta}))$ and
$(\theta^*,t^*(\theta^*))$ such that
\begin{eqnarray}\label{Taylor 1 E1}
  0 & = & P_n\frac{\partial}{\partial t}m\left(\theta^*,t^*(\theta^*)\right)+
  \left[\left(P_n\frac{\partial^2}{\partial t^2}m(\overline{\theta},\overline{c})\right)^{\top},
  \left(P_n\frac{\partial^2}{\partial \theta\partial t}
  m(\overline{\theta},\overline{c})\right)^{\top}\right]a_n\nonumber\\
\end{eqnarray}
with
\begin{equation}\label{a n }
    a_n:={\left({\left(\widehat{t}(\widehat{\theta})-t^*(\theta^*)\right)}^{\top},
{\left(\widehat{\theta}-\theta^*\right)}^{\top}\right)}^{\top}.
\end{equation}
\noindent  By Assumption \ref{assumption 7}, using  the UWLLN, we
can write
$$\left[P_n\frac{\partial^2}{\partial
t^2}m(\overline{\theta},\overline{c}),
  P_n\frac{\partial^2}{\partial\theta\partial t}
  m(\overline{\theta},\overline{c})\right]=\left[P_0\frac{\partial^2}{\partial t^2}
  m(\theta^*,t^*(\theta^*)),
 P_0\frac{\partial^2}{\partial\theta\partial t}
  m(\theta^*,t^*(\theta^*))\right]+o_P(1),$$
  to obtain
  from (\ref{Taylor 1 E1})
  \begin{equation}\label{Equation E1 donne}
    -P_n\frac{\partial}{\partial
    t}m(\theta^*,t^*)=
    \left[\left(P_0\frac{\partial^2}{\partial t^2}m(\theta^*,t^*)\right)^{\top}+o_P(1),
    \left(P_0\frac{\partial^2}{\partial\theta\partial
    t}m(\theta^*,t^*)\right)^{\top}
    +o_P(1)\right]a_n.
\end{equation}
In the same way, using a Taylor expansion in (\ref{E2}), we obtain
  \begin{equation}\label{Equation E2 donne}
    -P_n\frac{\partial}{\partial
    \theta}m(\theta^*,t^*)=
    \left[\left(P_0\frac{\partial^2}{\partial t\partial\theta}m(\theta^*,t^*)
    \right)^{\top}+o_P(1),
    \left(P_0\frac{\partial^2}{\partial
  \theta^2}m(\theta^*,t^*)\right)^\top+o_P(1)\right]a_n.
\end{equation}
From (\ref{Equation E1 donne}) and (\ref{Equation E2 donne}), we
get
\begin{eqnarray}\label{racine de n a n 1}
\sqrt{n}a_n  & = &  \sqrt{n}\left(
\begin{array}{cc}
  P_0\frac{\partial ^2}{\partial t^2}m(\theta^*,t^*) & \left(P_0\frac{\partial^2}
  {\partial\theta\partial t} m(\theta^*,t^*)\right)^{\top}\\
  \left(P_0\frac{\partial^2}{\partial t\partial \theta} m(\theta^*,t^*)\right)^{\top} &
  P_0\frac{\partial^2}{\partial \theta^2}m(\theta^*,t^*) \\
\end{array}
\right)^{-1}\times\nonumber\\
& & \times\left(
\begin{array}{c}
  -P_n\frac{\partial}{\partial t}m(\theta^*,t^*) \\
  -P_n\frac{\partial }{\partial \theta} m(\theta^*,t^*)\\
\end{array}
\right)+o_P(1).
\end{eqnarray}
Denote $S$ the $(1+l+d)\times(1+l+d)-$matrix defined by
\begin{equation}\label{S}
    S:=\left(
\begin{array}{cc}
  S_{11} & S_{12} \\
  S_{21} & S_{22} \\
\end{array}
\right):=\left(
\begin{array}{cc}
  P_0\frac{\partial ^2}{\partial t^2}m(\theta^*,t^*) & \left(P_0\frac{\partial^2}
  {\partial\theta\partial t} m(\theta^*,t^*)\right)^{\top}\\
  \left(P_0\frac{\partial^2}{\partial t\partial \theta} m(\theta^*,t^*)\right)^{\top} &
  P_0\frac{\partial^2}{\partial \theta^2}m(\theta^*,t^*) \\
\end{array}
\right).
\end{equation}
Hence, we obtain
\begin{equation*}
\sqrt{n}\left(
\begin{array}{c}
  \widehat{t}(\widehat{\theta})-t^*(\theta^*) \\
  \widehat{\theta} -\theta^*\\
\end{array}
\right)=\sqrt{n}S^{-1}\left(
\begin{array}{c}
  -P_n\frac{\partial}{\partial t}m(\theta^*,t^*) \\
  -P_n\frac{\partial}{\partial \theta}m(\theta^*,t^*)\\
\end{array}
\right)+o_P(1),
\end{equation*}
and the CL and Slutsky theorems conclude the proof.\\
 2) Using the fact that
 $$\widehat{t}(\widehat{\theta})-t^*(\theta^*)=O_P(1/\sqrt{n}),~
 P_n\partial m(\theta^*,t^*(\theta^*))/\partial t=
 P_0\partial m(\theta^*,t^*(\theta^*))/\partial t
 +o_P(1)=o_P(1)$$
 and  $$\widehat{\theta}-\theta^*=O_P(1/\sqrt{n}),~
 P_n\partial m(\theta^*,t^*(\theta^*))/\partial \theta=
 P_0\partial m(\theta^*,t^*(\theta^*))/\partial \theta
 +o_P(1)=o_P(1),$$
 we can write
 \begin{eqnarray}
  \sqrt{n}\left(\widehat{D}_\varphi(\mathcal{M},P_0)-D_\varphi(\mathcal{M},P_0)\right)
  & = &
  \sqrt{n}\left(P_nm(\widehat{\theta},\widehat{t}(\widehat{\theta}))-P_0m(\theta^*,t^*(\theta^*))\right)\nonumber\\
  & = & \sqrt{n}\left(P_nm(\theta^*,t^*(\theta^*))-P_0m(\theta^*,t^*(\theta^*))\right)+o_P(1),\nonumber
 \end{eqnarray}
 and the CL and Slutsky theorems end the proof.
\bibliographystyle{natbib}

\begin{thebibliography}{}
\bibitem[Baggerly(1998)]{Baggerly1998}
Baggerly, K.~A. (1998).
\newblock Empirical likelihood as a goodness-of-fit measure.
\newblock {\em Biometrika}, {\bf 85}(3), 535--547.

\bibitem[Bertail(2006)]{Bertail2006}
Bertail, P. (2006).
\newblock Empirical likelihood in some semiparametric models.
\newblock {\em Bernoulli}, {\bf 12}(2), 299--331.

\bibitem[Broniatowski and Keziou(2006)]{Bronia_Kez2006_STUDIA}
Broniatowski, M. and Keziou, A. (2006).
\newblock Minimization of {$\phi$}-divergences on sets of signed measures.
\newblock {\em Studia Sci. Math. Hungar.; arXiv:1003.5457}, {\bf 43}(4),
  403--442.

\bibitem[Broniatowski and Keziou(2008)]{BK_2006_div_emp}
Broniatowski, M. and Keziou, A. (2008).
\newblock Estimation and tests for models satisfying linear constraints with
  unknown parameter.
\newblock {\em arXiv:0811.3477v1}.

\bibitem[Broniatowski and Keziou(2009)]{Broniatowski_Keziou2009}
Broniatowski, M. and Keziou, A. (2009).
\newblock Parametric estimation and tests through divergences and the duality
  technique.
\newblock {\em J. Multivariate Anal.}, {\bf 100}(1), 16--36.

\bibitem[Chen {\em et~al.}(2007)]{Chen_Hong_Shum2007}
Chen, X., Hong, H., and Shum, M. (2007).
\newblock Nonparametric likihood ratio model selection tests between parametric
  likelihood and moment condition models.
\newblock {\em J. Econometrics}, {\bf 141}(1), 109--140.

\bibitem[Corcoran(1998)]{Corcoran1998}
Corcoran, S. (1998).
\newblock Bertlett adjustement of empirical discrepancy statistics.
\newblock {\em Biometrika}, {\bf 85}, 967--972.

\bibitem[Cressie and Read(1984)]{Cressie-Read1984}
Cressie, N. and Read, T. R.~C. (1984).
\newblock Multinomial goodness-of-fit tests.
\newblock {\em J. Roy. Statist. Soc. Ser. B}, {\bf 46}(3), 440--464.

\bibitem[Csisz{\'a}r(1963)]{Csiszar1963}
Csisz{\'a}r, I. (1963).
\newblock Eine informationstheoretische {U}ngleichung und ihre {A}nwendung auf
  den {B}eweis der {E}rgodizit\"at von {M}arkoffschen {K}etten.
\newblock {\em Magyar Tud. Akad. Mat. Kutat\'o Int. K\"ozl.}, {\bf 8}, 85--108.

\bibitem[Csisz{\'a}r(1967)]{Csiszar1967}
Csisz{\'a}r, I. (1967).
\newblock On topology properties of {$f$}-divergences.
\newblock {\em Studia Sci. Math. Hungar.}, {\bf 2}, 329--339.

\bibitem[Haberman(1984)]{Haberman1984}
Haberman, S.~J. (1984).
\newblock Adjustment by minimum discriminant information.
\newblock {\em Ann. Statist.}, {\bf 12}(3), 971--988.

\bibitem[Hansen {\em et~al.}(1996)]{Hansen_Healton_Yaron1996}
Hansen, L., Heaton, J., and Yaron, A. (1996).
\newblock Finite-sample properties of some alternative gmm estimators.
\newblock {\em Journal of Business and Economic Statistics}, {\bf 14},
  462--2800.

\bibitem[Hansen(1982)]{Hansen1982}
Hansen, L.~P. (1982).
\newblock Large sample properties of generalized method of moments estimators.
\newblock {\em Econometrica}, {\bf 50}(4), 1029--1054.

\bibitem[Hjort {\em et~al.}(2009)]{Hjort_McKeague_VanKeilegom2009}
Hjort, N.~L., McKeague, I.~W., and Van~Keilegom, I. (2009).
\newblock Extending the scope of empirical likelihood.
\newblock {\em Ann. Statist.}, {\bf 37}(3), 1079--1111.

\bibitem[Imbens(1997)]{Imbens1997}
Imbens, G.~W. (1997).
\newblock One-step estimators for over-identified generalized method of moments
  models.
\newblock {\em Rev. Econom. Stud.}, {\bf 64}(3), 359--383.

\bibitem[Keziou(2003)]{Keziou2003}
Keziou, A. (2003).
\newblock Dual representation of {$\phi$}-divergences and applications.
\newblock {\em C. R. Math. Acad. Sci. Paris}, {\bf 336}(10), 857--862.

\bibitem[Kitamura(2001)]{Kitamura2001}
Kitamura, Y. (2001).
\newblock Asymptotic optimality of empirical likelihood for testing moment
  restrictions.
\newblock {\em Econometrica}, {\bf 69}(6), 1661--1672.

\bibitem[Kitamura(2007)]{Kitamura2007}
Kitamura, Y. (2007).
\newblock {\em Empirical likelihood methods in econometric theory and
  practice}.
\newblock Cambridge University Press.

\bibitem[Liese and Vajda(1987)]{Liese-Vajda1987}
Liese, F. and Vajda, I. (1987).
\newblock {\em Convex statistical distances}, volume~95.
\newblock BSB B. G. Teubner Verlagsgesellschaft, Leipzig.

\bibitem[McCullagh and Nelder(1983)]{McCullagh_Nelder1983}
McCullagh, P. and Nelder, J.~A. (1983).
\newblock {\em Generalized linear models}.
\newblock Monographs on Statistics and Applied Probability. Chapman \& Hall,
  London.

\bibitem[Newey and Smith(2004)]{NeweySmith2004}
Newey, W.~K. and Smith, R.~J. (2004).
\newblock Higher order properties of {GMM} and generalized empirical likelihood
  estimators.
\newblock {\em Econometrica}, {\bf 72}(1), 219--255.

\bibitem[Owen(1990)]{Owen1990}
Owen, A. (1990).
\newblock Empirical likelihood ratio confidence regions.
\newblock {\em Ann. Statist.}, {\bf 18}(1), 90--120.

\bibitem[Owen(1991)]{Owen1991}
Owen, A. (1991).
\newblock Empirical likelihood for linear models.
\newblock {\em Ann. Statist.}, {\bf 19}(4), 1725--1747.

\bibitem[Owen(1988)]{Owen1988}
Owen, A.~B. (1988).
\newblock Empirical likelihood ratio confidence intervals for a single
  functional.
\newblock {\em Biometrika}, {\bf 75}(2), 237--249.

\bibitem[Owen(2001)]{Owen2001}
Owen, A.~B. (2001).
\newblock {\em Empirical Likelihood}.
\newblock Chapman and Hall, New York.

\bibitem[Pardo(2006)]{PardoLeandro2006}
Pardo, L. (2006).
\newblock {\em Statistical inference based on divergence measures}, volume 185
  of {\em Statistics: Textbooks and Monographs}.
\newblock Chapman \& Hall/CRC, Boca Raton, FL.

\bibitem[Qin and Lawless(1994)]{Qin-Lawless1994}
Qin, J. and Lawless, J. (1994).
\newblock Empirical likelihood and general estimating equations.
\newblock {\em Ann. Statist.}, {\bf 22}(1), 300--325.

\bibitem[Rockafellar(1970)]{Rockafellar1970}
Rockafellar, R.~T. (1970).
\newblock {\em Convex analysis}.
\newblock Princeton University Press, Princeton, N.J.

\bibitem[Schennach(2007)]{Schennach2007}
Schennach, S.~M. (2007).
\newblock Point estimation with exponentially tilted empirical likelihood.
\newblock {\em Ann. Statist.}, {\bf 35}(2), 634--672.

\bibitem[Sheehy(1987)]{Sheehy1987}
Sheehy, A. (1987).
\newblock \text{Kullback-Leibler} constrained estimation of probability
  measures.
\newblock {\em Report, Dept. Statistics, Stanford Univ.}

\bibitem[Smith(1997)]{Smith1997}
Smith, R.~J. (1997).
\newblock Alternative semi-parametric likelihood approches to generalized
  method of moments estimation.
\newblock {\em Economic Journal}, {\bf 107}, 503--519.

\bibitem[van~der Vaart(1998)]{vanderVaart1998}
van~der Vaart, A.~W. (1998).
\newblock {\em Asymptotic statistics}.
\newblock Cambridge Series in Statistical and Probabilistic Mathematics.
  Cambridge University Press, Cambridge.

\bibitem[White(1982)]{White1982}
White, H. (1982).
\newblock Maximum likelihood estimation of misspecified models.
\newblock {\em Econometrica}, {\bf 50}(1), 1--25.
\end{thebibliography}

\end{document}